\documentclass[letterpaper]{article}
\usepackage[dvips]{graphicx}
\usepackage{amsmath}

\begin{document}
{\large \bf On partitions of the unit interval generated by Brocot sequences. }\\

\medskip
\begin{center}Dushistova Anna \footnote{The work was supported by RFFI N06-01-00518-a. }\\
\end{center}
\medskip
 
Let $ p_{i, n}, i=1, \dots, 2^{n-1} $ be the lengths of intervals between the neighboring fractions of Brocot sequence $ F_n $. 
The asymptotic formula for 
 $ \sigma \left( F_n \right)=\sum_{i=1}^{N \left( n \right)}p_{i, n}^{\beta} $, 
improving known estimations, is obtained. \\

\medskip

{\large \bf\S1. Basic definitions and statements. }

\medskip

In the present work the partition of $ [0, 1] $ by the points of Brocot sequence is considered. \\
Brocot sequences $ F_n $, $ n=1, 2, \dots $ are defined inductively in the following way. \\
When $ n=1 $ let $ F_1=\{0, 1\}=\{\frac{0}{1}, \frac{1}{1}\} $. 
Let $ n\ge1 $ and for each $ k \le n $ sets $ F_k $ have been defined. 
Let's define $ F_{n+1} $. 
Consider fractions from $ F_n $, ordered by increase :
\begin{equation}
\label{ordered_elements}
0=x_{0, n}<x_{1, n}< \dots <x_{N \left( n \right), n}=1, N(n)=2^{n-1}.
\end{equation}
Then
$$ 
F_{n+1} = F_n \cup Q_{n+1}, 
$$ 
where $ Q_{n+1} $ is the set of mediants of neighboring fractions in $ F_n $, the given 
$$ 
Q_{n+1}=\{x_{i, n}\oplus x_{i-1, n}, i=1, \dots, N(n)\}, 
$$ 
where $ \frac{p}{q}\oplus \frac{p'}{q'}=\frac{p+p'}{q+q'} $. 
Elements in $ Q_n $ are known as Brocot fractions of order $ n $. Brocot sequences (known also as Stern-Brocot sequences) appeared in [1], [2]. Main properties of Brocot sequences can be found in [3],pages 140-143. \\ 
Let us consider the partition of $ [0, 1] $ by fractions of $ F_n $, that is with the given points like (\ref{ordered_elements}), 
let
 $ p_{i, n}=x_{i, n}-x_{i-1, n}, i=1, \dots, N \left( n \right) $ be lengths of $ [x_{i-1, n}, x_{i, n}) $. 
For fixed $ \beta $ we denote
$$ 
\sigma_{\beta} \left( F_n \right)
=\sum_{i=1}^{N \left( n \right)}p_{i, n}^{\beta}. 
$$ 
N. Moshchevitin and A.Zhigljavsky in [5] investigated the behavior of $ \sigma \left( F_n \right) $ 
when $ n $ tends to infinity. 
The following asymptotic equality was proved there.

{\bf Theorem 1. }\\
{\it For any $ \beta > 1 $ 
$$ 
\sigma_{\beta} \left( F_n \right) = \frac{2}{n^{\beta}}\frac{\zeta \left( 2\beta-1 \right)}{\zeta \left( 2\beta \right)}
+O \left( \frac{\log \left( n \right)}{n^{ \left( \beta+1 \right) \left( 2\beta-1 \right)/ \left( 2\beta \right)}} \right), n \to \infty, 
$$ 
where $ \zeta \left( s \right) $ is Riemann $ \zeta $ -function. 
}\\
The main result of this work is proof of the following more precise theorem.\\
{\bf Theorem 2. }\\
{\it For any $ \beta > 1 $ holds
$$ 
\sigma_{\beta} \left( F_n \right) =
\frac{1}{n^{\beta}}
\frac{2\zeta \left( 2\beta-1 \right)}{\zeta \left( 2\beta \right)}
+\sum_{1\le k < 2\beta -2}C_k\frac{1}{n^{\beta+k}}
+
\sum_{0\le k< \beta-2}{C^*}_k\frac{1}{n^{2\beta+k}}
+O \left( 
\frac{\log^{3\beta} n}{n^{3\beta-2}}
 \right), 
$$ 
where $ C_k \left( \beta \right), 1\le k \le 2\beta-2 $, $ {C^*}_k(\beta), 0 \le k\le \beta -2 $ are positive constants, 
depending on $ \beta $. 
}\\
When $ \beta \in ( 1, 1. 5] $ the formula in the theorem 2 is 
$$ 
\sigma_{\beta} \left( F_n \right) =
\frac{1}{n^{\beta}}
\frac{2\zeta \left( 2\beta-1 \right)}{\zeta \left( 2\beta \right)}
+O \left( 
\frac{\log^{3\beta} n}{n^{3\beta-2}}
 \right). 
$$ 
The error term here is better than in theorem 1, because when $ 1<\beta\leq 1. 5 $ 
 we have $ 3\beta-2 > \frac{ \left( \beta+1 \right) \left( 2\beta-1 \right)}{ \left( 2\beta \right)} $. 
\\When $ \beta > 1.5 $ theorem 2 gives the additional terms in asymptotic. \\
Note, that the history of the problem and review of some results is presented in the introduction of [5]. 

\medskip
\medskip
{\large \bf\S2. Some notation and formulation of auxiliary result.\\ }

\medskip

It's well known, that the sum of partial quotients in the continued fraction representation for Brocot fractions of order $ n $ equals $ n $, 
i. e.
$$ 
Q_n=\{\frac{p}{q}=[a_1, \dots, a_t], a_t \ge 2, a_1+ \dots +a_t=n \}. 
$$ 
Let $ A $ be the set of all integer vectors
 $ a= \left( a_1, \dots, a_t \right), t\ge1, a_j\ge1, j=1, \dots, t-1 $ and $ a_t\ge 2 $. \\
Let
$$ 
A_n=\{a= \left( a_1, \dots, a_t \right) \in A|a_1+ \dots +a_t=n\}. 
$$ 
Each $ a= \left( a_1, \dots, a_t \right) \in A $ is associated with the continued fraction 
 $ [0;a_1, \dots, a_t] $ (as integer part always equals zero, we simply denote it as $ [a_1, \dots, a_t] $ ) 
and corresponding continuant 
 $ \langle a_1, \dots, a_t\rangle $, empty continuant equals $ 1 $, -1 continuant 
equals 0. 
By construction, for any $ n>1 $ each fraction in $ F_n \backslash \left( F_1 \cup Q_n \right) $ 
has two neighbors in $ Q_n $, and each fraction $ \frac{p}{q} \in Q_n $ 
has two neighbors
 $ \frac{p_{-}}{q_{-}} $ and $ \frac{p_{+}}{q_{+}} $ in $ F_n \backslash Q_n $. 

{\bf Lemma 1. }
\newline
{\it For each $ a \in A_n $, the fraction $ \frac{p}{q} \in Q_n $ with denominator equal
to continuant $ q=\langle a_1, \dots, a_t \rangle $ has two neighbors
in $ F_n $ with denominators, equal to continuants $ q_{-}=\langle a_1, \dots, a_{t-1} \rangle $ and
 $ q_{+}=\langle a_1, \dots, a_{t}-1 \rangle $. 
Similarly, any fraction $ \frac{p}{q}\in F_{n-1}\backslash F_1 $ with denominator equal
to continuant
 $ \langle a_1, \dots, a_{t} \rangle $ has two neighbors in $ F_n $ with denominators, equal to continuants
$ \langle a_1, \dots, a_t, n- \left( a_1+ \dots +a_t \right)\rangle $ and
 $ \langle a_1, \dots, a_t-1, 1, n- \left( a_1+ \dots +a_t \right)\rangle $. }
\newline 
Proof is a simple induction with respect to $ n $ (see. [5]). 

To prove theorem 2 we need the following auxiliary result, 
that can be of self-contained interest. \\ 
Let
$$ 
\sigma_{\beta}(n)=\sum_{ \left( a_1, \dots, a_t \right) \in A_n}\frac{1}{\langle a_1, \dots, a_t\rangle^{2\beta}}
$$ 
with the fixed $ \beta > 1 $. 
\\
{\bf Theorem 3. }
\newline {\it For each $ \beta > 1 $ with some positive constants $ {C'}_k $, depending on $ \beta $, 
holds
$$ 
\sigma_{\beta}(n) 
=
\frac{1}{n^{2\beta}} \left( \frac{\zeta \left( 2\beta-1 \right)}{\zeta \left( 2\beta \right)}+2 \left( \frac{\zeta \left( 2\beta-1 \right)}{\zeta \left( 2\beta \right)} \right)^2 \right)+
\sum_{1\le k<2\beta-2}{C'}_k\frac{1}{n^{2\beta+k}}
+
O \left( 
\frac{\log^{4\beta} n}{n^{4\beta-2}}
 \right), 
$$ 
}
\\
In fact, in order to prove theorem 2 it is sufficient to obtain the main term in asymptotic in
theorem 3. This main term will be obtained in lemma 9 further. 
Note, that lemma 9 is the weaker variant of theorem 3 and it's actually used to prove theorem 3 in all completeness. \\

\medskip
\medskip
{\large \bf\S3. \bf Auxiliary statements. }\\
\medskip

{\bf Proof of theorem 3} uses 
splitting of $ \sigma $ which is the sum over $ A_n $ 
into the sums over smaller subsets of indices. \\
Let $ r $ and $ w $ be some integers, satisfying
the conditions $ r \ge 1 $ and $ 1 \le w \le n $. 
Then $ A_n=A_{n, 1}^{ \left( 1 \right)} \sqcup A_{n, 2}^{ \left( 1 \right)} $, where
$$ 
A_{n, 1}^{ \left( 1 \right)} = \{ a= \left( a_1, \ldots, a_t \right) \in A_n
 \mid \langle a_1, \ldots, a_t\rangle < n^r\}
$$ 
$$ 
A_{n, 2}^{ \left( 1 \right)} = A_n \backslash A_{n, 1}^{ \left( 1 \right)} = \{ a \in A_n
 \mid \langle a_1, \ldots, a_t\rangle \geq n^r\}
.$$

Then split $ A_{n, 1}^{ \left( 1 \right)} $ into 
$$ 
A_{n, 1}^{ \left( 2 \right)} = \{ a \in A_{n, 1}^{ \left( 1 \right)}
 \mid \max_{1 \leq j \leq t} a_j > n - w\}
$$ 
and
$$ 
A_{n, 2}^{ \left( 2 \right)} = A_{n, 1}^{ \left( 1 \right)} \backslash A_{n, 1}^{ \left( 2 \right)} = \{ a \in A_{n, 1}^{ \left( 1 \right)}
 \mid \max_{1 \leq j \leq t} a_j \leq n - w\}
.$$ 

Thus, all $ a \in A_{n, 1}^{ \left( 2 \right)} $ has at least one 
very large partial quotient; 
on the other hand, all $ a_j $ for $ a \in A_{n, 2}^{ \left( 2 \right)} $ are relatively small. 

So, 
 $ A_n=A_{n, 1}^{ \left( 2 \right)} \sqcup A_{n, 2}^{ \left( 2 \right)} \sqcup A_{n, 2}^{ \left( 1 \right)} $. 

For subset $ A_{n, i}^{ \left( j \right)} $ of $ A $ denote
$$ 
\Sigma_{n, i}^{ \left( j \right)} = \sum_{a \in A_{n, i}^{ \left( j \right)}}
 \frac{1} {q^{2\beta}}. 
$$ 

Thus, 
$$ 
\sigma_{\beta}(n) = 
\Sigma_{n, 2}^{ \left( 1 \right)} +
\Sigma_{n, 1}^{ \left( 2 \right)} +
\Sigma_{n, 2}^{ \left( 2 \right)}.
$$
Let us estimate these sums separately. 
\newline
{\bf Lemma 2. }
\newline
{\it Let $ n \geq 2, a = \left( a_1, \ldots, a_t \right) \in A_n $. 
 We have
 $ q=q_{+}+q_{-} \leq n q_{-}, $ 
 $ q_{-} \leq q_{+} \leq a_t q_{-}, $ 

$$ 
\sum_{a \in A_n} \left( \frac{1}{qq_{-}} + \frac{1}{qq_{+}} \right) = 1. 
$$ 
}
\\
The proof of this lemma is given in [3]. 
\\
\\
{\bf Lemma 3. }
\newline
{\it For each $ n \geq 1 $ holds
\begin{equation}
\label{large_cont}
\Sigma_{n, 2}^{ \left( 1 \right)} \leq
    2^{\beta-1}\frac{1}{n^{2r \left( \beta-1 \right)}}. 
\end{equation}
}
\\ ( Lemma 2 is similar to lemma 3 from [5]. )\\
{\bf Proof. }
\newline
As $ q \left( a \right)=\langle a_1, \dots, a_t\rangle \geq n^r $ for each $ a $ in $ A_{n, 2}^{ \left( 1 \right)} $, then 
using lemma 2, we obtain

$$ 
\Sigma_{n, 2}^{ \left( 1 \right)} 
\leq
\sum_{A^{ \left( 1 \right)}_{n, 2}}\frac{1}{ \left( qq_{+} \right)^{\beta}}
\leq
\max_{A^{ \left( 1 \right)}_{n, 2}}\frac{1}{ \left( qq_{+} \right)^{\beta-1}}\sum_{A^{ \left( 1 \right)}_{n, 2}}\frac{1}{qq_{+}}
\leq
$$ 
$$ 
\leq
\max_{A^{ \left( 1 \right)}_{n, 2}}\frac{1}{ \left( qq_{+} \right)^{\beta-1}}\sum_{A^{ \left( 1 \right)}_{n, 2}} \left( \frac{1}{qq_{+}}+\frac{1}{qq_{-}} \right)
\leq
\max_{A^{ \left( 1 \right)}_{n, 2}}\frac{1}{ \left( qq_{+} \right)^{\beta-1}}
\leq
$$ 
$$ 
\leq
\frac{2^{\beta-1}}{q^{2 \left( \beta-1 \right)}}
\leq
\frac{2^{\beta-1}}{n^{2r \left( \beta-1 \right)}}.
$$ 
{\bf Lemma is proved. }
\newline
{\bf Lemma 4. }
\newline
{\it For each $ a \in A_{n, 1}^{ \left( 1 \right)} $ when $ n \geq 2 $ holds

$$ 
t \leq Kr \log n, 
 K = \left( \log \frac{\sqrt{5}+1}{2} \right)^{-1}. 
$$ }

{\bf Proof. }
\newline
For each $ a \in A_{n, 1}^{ \left( 1 \right)} $ holds
$$ 
 \left( \frac{\sqrt{5}+1}{2} \right)^{t}
\leq
\langle a_1, \ldots, a_t\rangle
\leq
n^r
, \quad
t \left( \log \frac{\sqrt{5}+1}{2} \right)
\leq
r \log n . 
$$ 
{\bf Lemma is proved. }\\
\newline
{\bf Lemma 5. }\\
{\it
When $ n \rightarrow \infty $, the following estimate for $ \sum_{n, 2}^{ \left( 2 \right)} $ holds:
\begin{equation}
\label{small_elem}
\Sigma_{n, 2}^{ \left( 2 \right)}
\ll
\frac{n^2\log^{4\beta}n}{w^{4\beta}}. 
\end{equation}
}
Note, that lemma 5 can be improved, 
but it won't effect the main result. 
{\bf Proof. }
\newline
According to lemma 4 for each $ a \in A_{n, 1}^{ \left( 1 \right)} $ holds $ t \le K r \log n $. 
As $ n = a_1 + \ldots + a_t \leq t \max a_j $, then
 $ \max a_j \geq \frac{n}{K r \log n} $. 
\newline
Let $ a \in A_{n, 2}^{ \left( 2 \right)} $ and $ j $ be such that
 $ a_j = \max \{ a_1, \ldots, a_t\}. $ 
As $ a_j \leq n - w $, then for the sum of the rest $ a_j $ we have 
 $ \sum\limits_{i \not= j} a_j \ge w $, and, similarly to the above, 
 $ \max_{i \not= j} a_j \geq \frac{w}{K r \log n} $. 
\newline
This implies, that for each $ a \in A_{n, 2}^{ \left( 2 \right)} $ there exist 2 different
partial quotients $ a_k $ and $ a_l, k \neq l $, that $ a_k \geq \frac{w}{K r \log n} $, $ a_l \geq \frac{w}{K r \log n} $. 
Hence, 

$$ 
\Sigma_{n, 2}^{ \left( 2 \right)} \le
\sum_{
\begin{array}{c}
a_1+\ldots+a_t = n, \\
\langle a_1, \dots, a_t\rangle \leq n^r, \\
\exists k, l, k\ne l : a_k, a_l \geq \frac{w}{K\log n}
\end{array}
}
\frac{1}{q^{2\beta}}. 
$$ 
Using the well known formula for continuants (see. [3]), we get
\begin{equation}
\label{knut}
\begin{split}
\langle a_1, \dots, a_i, \dots, a_t\rangle=
a_i\langle a_1, \dots, a_{i-1}\rangle\langle a_{i+1}, \dots, a_t\rangle
\times
\\
\times
 \left( 1+\frac{1}{a_i}[a_{i-1}, \dots, a_1]+ \frac{1}{a_i}[a_{i+1}, \dots, a_t] \right)=
%$$  
\\
%$$  
=a_i\langle a_1, \dots, a_{i-1}\rangle\langle a_{i+1}, \dots, a_t\rangle
 \left( 1+\frac{1}{a_i}[a_{i-1}, \dots, a_1] \right)
\times
\\
\times
 \left( 1+[a_{i+1}, \dots a_t]\frac{1}{a_i+[a_{i-1}, \dots, a_1]} \right), 
\end{split}
\end{equation}
Therefore, 
$$ 
q \left( a \right) \geq a_ka_l\langle a_1, . . , a_{\min \left( k, l \right)-1}\rangle
\langle a_{\min \left( k, l \right)+1}, . . , a_{\max \left( k, l \right)-1}\rangle
\langle a_{\max \left( k, l \right)+1}, . . , a_t\rangle. 
$$ 
Note, that with the $ X $ fixed the elements of set 
$$ 
\{a \in A_n |\langle a_1, \dots, a_t\rangle \leq n^r, 
\exists k, l, k \neq l : a_k, a_l \geq X \}
$$ 
look like
$$
\left( a_1, \dots a_{\min \left( k, l \right)-1}, T, a_{\min \left( k, l \right)+1}, \dots,
a_{\max \left( k, l \right)-1}, P, a_{\max \left( k, l \right)+1}, \dots, a_t \right),
$$
where $ T, P \geq X $, lengths of 
$$
\left( a_1, \dots, a_{\min \left( k, l \right)-1} \right),
$$
$$
\left( a_{\min \left( k, l \right)+1}, \dots, a_{\max \left( k, l \right)-1} \right)
$$
and
$$ 
\left( a_{\max \left( k, l \right)+1}, \dots, a_t \right) 
$$
are not fixed, and the sum is
$$
a_1+ \dots +a_{l-1}+a_{l+1}+ \dots +a_{k-1}+a_{k+1}+ \dots +a_t=n-T-P.
$$
Let 
$$
a_1+ \dots +a_{\min \left( k, l \right)-1}=u,
$$
$$
a_{\min \left( k, l \right)+1}+ \dots +a_{\max \left( k, l \right)-1}=v,
$$
$$
a_{\max \left( k, l \right)+1}+ \dots +a_t=s,
$$
$$
u+v+s=n-T-P.
$$
Thus, 
$$ 
\Sigma_{n, 2}^{ \left( 2 \right)}
\leq
\sum_{
\begin{array}{c}
a_1+\ldots+a_t = n, \\
\langle a_1, \dots, a_t\rangle \leq n^r, \\
\exists k, l, k\ne l : a_k, a_l \geq \frac{w}{K\log n}
\end{array}
}
\frac{1}{ (a_k a_l)^{2\beta} }
\times
$$
$$
\times
\frac{1}{ \left( \langle a_1, . . , a_{\min \left( k, l \right)-1}\rangle
\langle a_{\min \left( k, l \right)+1}, . . , a_{\max \left( k, l \right)-1}\rangle
\langle a_{\max \left( k, l \right)+1}, . . , a_t\rangle \right)^{2\beta}}
\leq
$$ 
$$ 
\leq
\sum_{
\begin{array}{c}
a_k+a_l \le n, \\
a_k, a_l \geq \frac{w}{K\log n}
\end{array}}
\frac{1}{ \left( a_ka_l \right)^{2\beta}}
\times
$$ 
$$ 
\times
\sum_{u+v+s=n-a_k-a_l}
\sum_{a_1+ \dots +a_{\min \left( k, l \right)-1}=u}\frac{1}{\langle a_1, \dots, a_{\min \left( k, l \right)-1} \rangle^{2\beta}}
\times
$$ 
$$ 
\times
\sum_{a_{\min \left( k, l \right)+1}+ \dots +a_{\max \left( k, l \right)-1}=v}\frac{1}{\langle a_{\min \left( k, l \right)+1}, \dots, a_{\max \left( k, l \right)-1} \rangle^{2\beta}}
\times
$$ 
$$ 
\times
\sum_{a_{\max \left( k, l \right)+1}+ \dots +a_t=s}\frac{1}{\langle a_{\max \left( k, l \right)+1}, \dots, a_t\rangle^{2\beta}}.
$$ 
Let's estimate the internal sum. 
$$ 
\sum_{u+v+s=n-a_k-a_l}
\sum_{a_1+ \dots +a_{\min \left( k, l \right)-1}=u}\frac{1}{\langle a_1, \dots, a_{\min \left( k, l \right)-1} \rangle^{2\beta}}
$$ 
$$ 
\sum_{a_{\min \left( k, l \right)+1}+ \dots +a_{\max \left( k, l \right)-1}=v}\frac{1}{\langle a_{\min \left( k, l \right)+1}, \dots, a_{\max \left( k, l \right)-1} \rangle^{2\beta}}
\sum_{a_{\max \left( k, l \right)+1}+ \dots +a_t=s}\frac{1}{\langle a_{\max \left( k, l \right)+1}, \dots, a_t\rangle^{2\beta}}=
$$ 
$$ 
=\sum_{u+v+s=n-a_k-a_l}
\sum_{x_1+ \dots +x_r=u}\frac{1}{\langle x_1, \dots, x_r \rangle^{2\beta}}
\sum_{y_1+ \dots +y_h=v}\frac{1}{\langle y_1, \dots, y_h \rangle^{2\beta}}
\sum_{z_1+ \dots +z_g=s}\frac{1}{\langle z_1, \dots, z_g\rangle^{2\beta}} \le
$$ 
$$ 
\le
\sum_{u+v+s \le n}
\sum_{x_1+ \dots +x_r=u}\frac{1}{\langle x_1, \dots, x_r \rangle^{2\beta}}
\sum_{y_1+ \dots +y_h=v}\frac{1}{\langle y_1, \dots, y_h \rangle^{2\beta}}
\sum_{z_1+ \dots +z_g=s}\frac{1}{\langle z_1, \dots, z_g\rangle^{2\beta}} \le
$$ 
$$ 
\le
\sum_{u+v+s \le \infty}
\sum_{x_1+ \dots +x_r=u}\frac{1}{\langle x_1, \dots, x_r \rangle^{2\beta}}
\sum_{y_1+ \dots +y_h=v}\frac{1}{\langle y_1, \dots, y_h \rangle^{2\beta}}
\sum_{z_1+ \dots +z_g=s}\frac{1}{\langle z_1, \dots, z_g\rangle^{2\beta}} =
$$ 
$$ 
=
\sum_{x_1+ \dots +x_r \le \infty}\frac{1}{\langle x_1, \dots, x_r \rangle^{2\beta}}
\sum_{y_1+ \dots +y_h \le \infty}\frac{1}{\langle y_1, \dots, y_h \rangle^{2\beta}}
\sum_{z_1+ \dots +z_g \le \infty}\frac{1}{\langle z_1, \dots, z_g\rangle^{2\beta}} \le
$$ 
$$ 
\le
 \left( \sum_{x_1+ \dots +x_r \le \infty}\frac{1}{\langle x_1, \dots, x_r \rangle^{2\beta}} \right)^3
 \le
 \left(1+ 2\sum_{
\begin{array}{c}
x_1+ \dots +x_r \le \infty, \\
 x_r \ge 2
\end{array}
}
\frac{1}{\langle x_1, \dots, x_r \rangle^{2\beta}} \right)^3=
$$ 
$$ 
= 
\left(1+2 \frac{\zeta \left( 2\beta-1 \right)}{\zeta \left( 2\beta \right)} \right)^3. 
$$ 
Then for the external sum the following estimation holds:
$$ 
\sum_{
\begin{array}{c}
a_k+a_l \le n, \\
a_k, a_l \geq \frac{w}{K\log n}
\end{array}
}
\frac{1}{ \left( a_ka_l \right)^{2\beta}} \ll
 \left( n-2\frac{w}{K \log n} \right)^2 \left( \frac{K^{4\beta}\log^{4\beta}n}{w^{4\beta}} \right), 
$$ 
where $ \left( n-2\frac{w}{K \log n} \right)^2 $ is the number of items in the sum, 
and $ \frac{K^{4\beta}\log^{4\beta}n}{w^{4\beta}} $ is the estimate for
 $ 
\frac{1}{ \left( a_ka_l \right)^{2\beta}}
 $. 
Thus, we obtain that
$$ 
\Sigma_{n, 2}^{ \left( 2 \right)} \ll
 \left( n-2\frac{w}{K \log n} \right)^2 \left( \frac{\log^{4\beta}n}{w^{4\beta}} \right)
=O \left( \frac{n^2\log^{4\beta}n}{w^{4\beta}} \right). 
$$ 

{\bf Lemma is proved. }
\\
{\bf Lemma 6. }
\\
{\it When $ w < \frac{n}{2} $ 

$$ 
Q=\{a\in A_n|\exists j: a_j>n-w\}=\bigsqcup_{X=n-w}^{n}\bigsqcup_{u+v=n-X}P \left( u, v, X \right), 
$$ 
where
$$ 
P \left( u, v, X \right)=\{a \in A|a= \left( a_1, \dots, a_t, X, {a'}_1, \dots, {a'}_{t'} \right),
$$ 
$$ 
a_1+ \dots +a_t=u, {a'}_1+ \dots +{a'}_{t'}=v\}, 
$$ 
and symbol $ \bigsqcup $ means, that sets $ P \left( u, v, X \right) $ and $ P \left( u', v', X' \right) $ don't intersect,
when
 $ \left( u, v, X \right) \ne \left( u', v', X' \right) $. 
}

{\bf Proof. }
\\
Let $ a \in Q $, then there exists the partial quotient $ a_j=Y>n-w $, hence
this element is in $ P \left( a_1+ \dots +a_{j-1}, a_{j+1}+ \dots +a_{t}, Y \right) $. \\
Inversely, if $ a \in \bigsqcup_{X=w}^{n}\bigsqcup_{u+v=n-X}P \left( u, v, X \right) $, then $ a \in A_n $ and
there exists the partial quotient $ a_i>n-w $. \\
 Let's prove that element from $ Q $ can't belong to several sets
 $ P \left( u, v, X \right) $
 at the same time.
 If it's not true and there exists $ a \in Q $, such that $ a \in P \left( u, v, X \right) $ and
 $ a \in P \left( u^{*}, v^{*}, X^{*} \right) $, then it can be represented as 
$$ 
a= \left( a_1, \dots, a_{i-1}, X, a_{i+1}, \dots , a_t \right), 
$$ 
$$ 
a_1+ \dots +a_{i-1}=u, 
a_{i+1}+ \dots +a_t=v, 
u+v=n-X
$$ 
and 
$$ 
a= \left( a^*_1, \dots , a^*_{j-1}, X^*, a^*_{j+1}, \dots , a^*_t \right), 
$$ 
$$ 
a^*_1+ \dots +a^*_{j-1}=u^*, 
a^*_{j+1}+ \dots +a^*_t=v^*, 
u^*+v^*=n-X^*. 
$$ 
Let $ i \ne j $. Then in $ a $ there exist two partial quotients, larger than $ \frac{n}{2} $, 
and hence $ \sum a_i > n $, that contradicts the fact that $ a \in Q $. \\
Hence, $ i = j $, i. e. $ X = X^{*} $, and, obviously, $ \left( u, v \right) = \left( u^{*}, v^{*} \right) $, the given sets $ P \left( u, v, X \right) $ and $ P \left( u^{*}, v^{*}, X^{*} \right) $ are the same. \\
{\bf Lemma is proved. }\\
{\bf Lemma 7. }
{\it
Let $ w < \frac{n}{2} $. \\
Then for the sum 
$$ 
R_0=
\sum_{
\begin{array}{c}
a \in A_n, \\
\exists j:a_j>n-w
\end{array}
}
\frac{1}{\langle a_1, \dots , a_{j-1}\rangle^{2\beta}\langle a_{j+1}, \dots , a_t\rangle^{2\beta}}
$$ 
the following asymptotic formula holds:
$$ 
R_0=
\sum_{
\begin{array}{c}
a \in A_n, \\
\exists j:a_j>n-w
\end{array}
}
\frac{1}{\langle a_1, \dots , a_{j-1}\rangle^{2\beta}\langle a_{j+1}, \dots , a_t\rangle^{2\beta}}
=C_0+O \left( \frac{1}{w^{2 \left( \beta-1 \right)}} \right), 
$$ 
where 
$$ 
C_0=\frac{\zeta \left( 2\beta-1 \right)}{\zeta \left( 2\beta \right)}+2 \left( \frac{\zeta \left( 2\beta-1 \right)}{\zeta \left( 2\beta \right)} \right)^2.
$$ 
}\\
{\bf Proof. }
\newline
According to lemma 6, 
$$ 
R_0=
\sum_{
\begin{array}{c}
a \in A_n, \\
\exists j:a_j>n-w
\end{array}
}
\frac{1}{\langle a_1, \dots , a_{j-1}\rangle^{2\beta}\langle a_{j+1}, \dots , a_t\rangle^{2\beta}}=
$$ 
$$ 
=\sum_{X=n-w}^{n}
\sum_{u+v=n-X}
\sum_{
\begin{array}{c}
a_1+ \dots +a_{j-1}=u, \\
a_{j+1}+ \dots +a_t=v
\end{array}
}
\frac{1}{\langle a_1, \dots , a_{j-1}\rangle^{2\beta}\langle a_{j+1}, \dots , a_t\rangle^{2\beta}}=
$$ 
$$ 
\sum_{u+v \le w}
\sum_{a_1+ \dots +a_{j-1}=u}
\frac{1}{\langle a_1, \dots , a_{j-1}\rangle^{2\beta}}
\sum_{
\begin{array}{c}
a_{j+1}+ \dots +a_t=v, \\
a_t\geq 2
\end{array}
}
\frac{1}{\langle a_{j+1}, \dots , a_t\rangle^{2\beta}}. 
$$ 
Let's split the sum into 2 parts, separating item $ \sum_1 $ with $ u=1 $. Replacing sum with $ u>1 $ to
the doubled $ \sum_2 $ over the set of indices with the last partial quotient larger or equal to 2, 
i.e. over $ A_n $, we obtain 
$$ 
R_0=\Sigma_1+2\Sigma_2, 
$$ 
where
$$ 
\Sigma_1=
\sum_{v \le w-1} 
\sum_{
\begin{array}{c}
a_{j+1}+ \dots +a_t=v, \\
a_t \geq 2
\end{array}}
\frac{1}{\langle a_{j+1}, \dots , a_t\rangle^{2\beta}}
=
\sum_{v \le w-1} 
\sigma_{\beta}(v)
, 
$$ 
$$ 
\Sigma_2=
\sum_{u+v \le w}
\sum_{
\begin{array}{c}
a_1+ \dots +a_{j-1}=u, \\
a_{j-1} \ge 2
\end{array}
}
\frac{1}{\langle a_1, \dots , a_{j-1}\rangle^{2\beta}}
\times
$$ 
$$ 
\times
\sum_{
\begin{array}{c}
a_{j+1}+ \dots +a_t=v, \\
a_t \ge 2
\end{array}
}
\frac{1}{\langle a_{j+1}, \dots , a_t\rangle^{2\beta}}=
$$ 
$$ 
=\sum_{u+v \le w}
\sigma_{\beta}(u)
\sigma_{\beta}(v). 
$$ 
For $ \Sigma_1 $, doing the calculations like in the proof of 
lemma 7 of the theorem 1 (see [3]), we get
$$ 
\Sigma_1=\frac{\zeta \left( 2\beta-1 \right)}{\zeta \left( 2\beta \right)}+O \left( \frac{1}{w^{2\beta-2}} \right). 
$$ 
Let's get the asymptotic formula for $ \Sigma_2 $ :
$$ 
\Sigma_2 \leq 
\sum_{
\begin{array}{c}
a_1+. . +a_{k} \leq \infty, \\
a_k \ge 2
\end{array}
}
\frac{1}{\langle a_1, \dots , a_k\rangle^{2\beta}}
\sum_{
\begin{array}{c}
a_1+. . +a_l \leq \infty, \\
a_l \ge 2
\end{array}
}
\frac{1}{\langle a_1, \dots , a_l\rangle^{2\beta}}=
$$ 
$$ 
=
 \left( \frac{\zeta \left( 2\beta-1 \right)}{\zeta \left( 2\beta \right)} \right) ^2. 
$$ 
On the other hand, 
$$ 
\Sigma_2 \geq
\sum_{
\begin{array}{c}
a_1+. . +a_{k} \leq \frac{w}{2}, \\
a_k \ge 2
\end{array}
}
\frac{1}{\langle a_1, \dots , a_k\rangle^{2\beta}}
\sum_{
\begin{array}{c}
a_1+. . +a_l \leq \frac{w}{2}, \\
a_l \ge 2
\end{array}
}
\frac{1}{\langle a_1, \dots , a_l\rangle^{2\beta}}, 
$$ 
because all these items are in $ \sum_2 $. 
$$ 
\Sigma_2 \geq
\sum_{
\begin{array}{c}
a_1+. . +a_{k} \leq \infty, \\
a_k \ge 2
\end{array}
}
\frac{1}{\langle a_1, \dots , a_k\rangle^{2\beta}}
\sum_{
\begin{array}{c}
a_1+. . +a_l \leq \infty, \\
a_l \ge 2
\end{array}
}
\frac{1}{\langle a_1, \dots , a_l\rangle^{2\beta}}
-
$$
$$
-
O \left( \frac{1}{w^{2\beta -2}} \right)=
$$ 
$$ 
= \left( \frac{\zeta \left( 2\beta-1 \right)}{\zeta \left( 2\beta \right)} \right) ^2+O \left( \frac{1}{w^{2 \left( \beta-1 \right)}} \right),
$$ 
i.e. 
$$ 
\Sigma_2 = \left( \frac{\zeta \left( 2\beta-1 \right)}{\zeta \left( 2\beta \right)} \right)^2+O \left( \frac{1}{w^{2 \left( \beta-1 \right)}} \right). 
$$ 
Thus, 
$$ 
R_0= C_0+O \left( \frac{1}{w^{2 \left( \beta-1 \right)}} \right). 
$$ 
{\bf Lemma is proved. }\\
{\bf Lemma 8. }
{\it When $ w \le \frac{n}{2}-2 $ 
\begin{equation}
\label{main_term_weak}
\Sigma_{n, 1}^{ \left( 2 \right)}
=\frac{C_0}{n^{2\beta}}
+O \left( \frac{w}{n^{2\beta+1}} \right)+O \left( \frac{1}{n^{2\beta}w^{2 \left( \beta-1 \right)}} \right)
+O \left( \frac{1}{n^{2r(\beta-1)}} \right). 
\end{equation}
where $C_0$ is defined in Lemma 7.
}\\

{\bf Proof. }
\\
Note that when $ w \le \frac{n}{2}-2 $, each element of $ A_{n, 1}^{ \left( 2 \right)} $ has the only one partial quotient, that is larger, than $ n-w $. 
$$ 
\Sigma_{n, 1}^{ \left( 2 \right)}
=
\sum_{
\begin{array}{c}
a \in A_n, \\
q \left( a \right) < n^r, \\
\exists i: a_i > n-w
\end{array}
}
\frac{1}{q^{2\beta}}=
$$ 
$$ 
=
\sum_{
\begin{array}{c}
a \in A_n, \\
\exists i: a_i > n-w
\end{array}
}
\frac{1}{q^{2\beta}}
-
\sum_{
\begin{array}{c}
a \in A_n, \\
q \left( a \right) \geq n^r, \\
\exists i: a_i > n-w
\end{array}
}
\frac{1}{q^{2\beta}}.
$$ 

The second sum is estimated according to Lemma 3. 
$$ 
\sum_{
\begin{array}{c}
a \in A_n, \\
q \left( a \right) \geq n^r, \\
\exists i: a_i > n-w
\end{array}
}
\frac{1}{q^{2\beta}}=
O \left( 
\frac{1}{n^{2r \left( \beta-1 \right)}}
 \right) 
$$ 
Let's estimate the first sum. 
Let $ a_i=n+O \left( w \right) $. \\
Using the formula for continuants (\ref{knut}), 
obtain
$$ 
q \left( a \right) = \left( a_i+[a_{i-1}, \dots , a_1]+[a_{i+1}, \dots , a_t] \right)\langle a_1, \dots , a_{i-1}\rangle \langle a_{i+1}, \dots , a_t\rangle=
$$ 
$$ 
=n \left( 1+\frac{w}{n}\cdot\theta \right)\langle a_1, \dots , a_{i-1}\rangle\langle a_{i+1}, \dots , a_t\rangle, 
$$ 
where $ \theta=\theta \left( w, n \right), |\theta|\le 1 $. 
Then, 
considering $ \frac{1}{q \left( a \right)^{2\beta}} $ as the function of argument $ a_i $ 
and expanding it in Taylor series according to argument $ \frac{w}{n}\cdot\theta $, 
we obtain 
$$ 
\frac{1}{q \left( a \right)^{2\beta}}=
\frac{1}{n^{2\beta}\langle a_1, \dots , a_{j-1}\rangle ^{2\beta} \langle a_{j+1}, \dots , a_t\rangle^{2\beta}}
 \left( 1+O \left( \frac{w}{n} \right) \right). 
$$ 

Thus, 
$$ 
\sum_{
\begin{array}{c}
a \in A_n, \\
\exists j: a_j > n-w
\end{array}
}
\frac{1}{q^{2\beta}}
= 
$$
$$
=
\sum_{
\begin{array}{c}
a \in A_n, \\
\exists j: a_j > n-w
\end{array}
}
\frac{1}{n^{2\beta}\langle a_1, \dots , a_{j-1}\rangle^{2\beta}\langle a_{j+1}, \dots , a_t \rangle^{2\beta}} \left( 1+O \left( \frac{w}{n} \right) \right)=
$$ 
$$ 
=\frac{1}{n^{2\beta}}
\sum_{
\begin{array}{c}
a \in A_n, \\
\exists j: a_j > n-w
\end{array}
}
\frac{1}{\langle a_1, \dots , a_{j-1}\rangle^{2\beta}\langle a_{j+1}, \dots , a_t \rangle^{2\beta}} \left( 1+O \left( \frac{w}{n} \right) \right). 
$$ 
Hence, according to lemma 7 
$$ 
\sum_{
\begin{array}{c}
a \in A_n, \\
\exists i: a_i > n-w
\end{array}
}
\frac{1}{q^{2\beta}}
=
\frac{1}{n^{2\beta}}
 \left( C_0+O \left( \frac{1}{w^{2 \left( \beta-1 \right)}} \right) \right)
 \left( 1+O \left( \frac{w}{n} \right) \right)=
$$ 
$$ 
=
\frac{C_0}{n^{2\beta}}
+O \left( \frac{w}{n^{2\beta+1}} \right)
+O \left( \frac{1}{n^{2\beta}w^{2 \left( \beta-1 \right)}} \right). 
$$ 
{\bf Lemma is proved. }\\
The following lemma is the weaker variant of theorem 3, which will be 
used to prove more precise result. \\
{\bf Lemma 9. }
{\it
When $ \beta > 1 $ we get
$$ 
\sigma_{\beta} \left( n \right)
= \frac{C_0}{n^{2\beta}}
+O \left( \frac{\log^{\frac{4\beta}{4\beta+1}} n}{n^{2\beta+1-\frac{2\beta+3}{4\beta+1}}} \right), 
$$ 
}

{\bf Proof. }\\
Using (\ref{large_cont}), (\ref{small_elem}), (\ref{main_term_weak}), we get
$$ 
\sigma_{\beta} \left( n \right) = 
\Sigma_{n, 2}^{ \left( 1 \right)} +
\Sigma_{n, 1}^{ \left( 2 \right)} +
\Sigma_{n, 2}^{ \left( 2 \right)}=\frac{C_0}{n^{2\beta}}
+O \left( \frac{1}{n^{2\beta}w^{2 \left( \beta-1 \right)}} \right)
+O \left( \frac{w}{n^{2\beta+1}} \right)+
$$ 
$$ 
+O \left( \frac{1}{n^{2r \left( \beta-1 \right)}} \right)+
O \left( \frac{n^2\log^{4\beta}n}{w^{4\beta}} \right). 
$$

Optimizing according to $ w $ and $ r $
$$
w =\min \left\{ \frac{n}{2}-2,
n^{\frac{2\beta+3}{4\beta+1}}\log^{\frac{4\beta}{4\beta+1}} n \right\},
$$
$$
r=\frac{2\beta+1}{2 \left( \beta-1 \right)},
$$
we obtain
$$
\sigma \left( F_n \right) = \frac{C_0}{n^{2\beta}}
+O \left( \frac{\log^{\frac{4\beta}{4\beta+1}} n}{n^{2\beta+1-\frac{2\beta+3}{4\beta+1}}} \right).
$$

{\bf Lemma is proved. }\\

\medskip
\medskip
{\large \bf\S4. Main lemma and final step of proving theorem 3. }
\medskip
\\
Let's consider the sum $ \Sigma_{n, 1}^{ \left( 2 \right)} $. 

{\bf Lemma 10. }
{\it Let $ w < \frac{n}{2} $. \\
In case $ 2\beta $ is integer, 
\begin{equation}
\label{main_term_asim_log}
\begin{split}
\Sigma_{n, 1}^{ \left( 2 \right)}
=
\frac{C_0}{n^{2\beta}}+{}\\
+
\sum_{1\le k<2\beta-2}{C'}_k\frac{1}{n^{2\beta+k}}
+
O \left( 
\frac{1}{n^{2\beta}w^{2 \left( \beta-1 \right)}}
+
\frac{\log n}{n^{4\beta-2}} 
+
\frac{1}{n^{2r \left( \beta-1 \right)}}
 \right), 
\end{split}
\end{equation}
in other case, 
\begin{equation}
\label{main_term_asim_nolog}
\begin{split}
\Sigma_{n, 1}^{ \left( 2 \right)}
=
\frac{C_0}{n^{2\beta}} +{}\\
+
\sum_{1\le k<2\beta-2}{C'}_k\frac{1}{n^{2\beta+k}}
+
O \left( 
\frac{1}{n^{2\beta}w^{2 \left( \beta-1 \right)}}
+
\frac{1}{n^{2r \left( \beta-1 \right)}}
 \right), 
\end{split}
\end{equation}
where $ {C'}_k $ are some constants. 
}\\

{\bf Final step of proving theorem 3. }
Using (\ref{large_cont}), (\ref{small_elem}), (\ref{main_term_asim_log}) and (\ref{main_term_asim_nolog}) 
error term $ R $ in case $ 2\beta $ is integer equals 
$$ 
R= O \left( 
\frac{1}{n^{2\beta}w^{2 \left( \beta-1 \right)}}
+
\frac{\log n}{n^{4\beta-2}}
+
\frac{1}{n^{2r \left( \beta-1 \right)}}
+
\frac{n^2 \log^{4\beta}n}{w^{4\beta}}
 \right), 
$$ 
in case $ 2\beta $ is not integer equals 
$$ 
R= O \left( 
\frac{1}{n^{2\beta}w^{2 \left( \beta-1 \right)}}
+
\frac{1}{n^{2r \left( \beta-1 \right)}}
+
\frac{n^2 \log^{4\beta}n}{w^{4\beta}}
 \right). 
$$ 
Let $ w =\frac{n}{2}-2 $, $ r=\frac{2\beta-1}{\beta-1} $. 
Substituting the value of $ w $ and $ r $, we get 
$$ 
R=O \left( 
\frac{\log^{4\beta} n}{n^{4\beta-2}}
 \right), 
$$ 
so theorem 3 follows. \\
{\bf Proof of lemma 10. }
\\
\begin{equation}
\label{1st}
\begin{split}
\Sigma_{n, 1}^{ \left( 2 \right)}
=
\sum_{
\begin{array}{c}
a \in A_n, \\
q \left( a \right) < n^r, \\
\exists i: a_i > n-w
\end{array}
}
\frac{1}{q^{2\beta}}
=
\sum_{
\begin{array}{c}
a \in A_n, \\
\exists i: a_i > n-w
\end{array}
}
\frac{1}{q^{2\beta}}
-
\sum_{
\begin{array}{c}
a \in A_n, \\
q \left( a \right) \geq n^r, \\
\exists i: a_i > n-w
\end{array}
}
\frac{1}{q^{2\beta}}
\end{split}
\end{equation}
Second sum in (\ref{1st}) can be estimated according to lemma 3 as
$$ 
\sum_{
\begin{array}{c}
a \in A_n, \\
q \left( a \right) \geq n^r, \\
\exists i: a_i > n-w
\end{array}
}
\frac{1}{q^{2\beta}}=
O \left( 
\frac{1}{n^{2r \left( \beta-1 \right)}}
 \right). 
$$ 

Let $ a_i=n-v $, where $ v=1, \dots , w-1 $. 
Using (\ref{knut}), 
we obtain
$$ 
q \left( a \right)=\langle a_1, \dots , a_{i-1}\rangle\langle a_{i+1}, \dots , a_t\rangle
 \left( a_i+[a_{i-1}, \dots , a_1]+[a_{i+1}, \dots a_t] \right)=
$$ 
$$ 
=n\langle a_1, \dots , a_{i-1}\rangle\langle a_{i+1}, \dots , a_t\rangle
 \left( 1-\frac{v}{n}+\frac{1}{n} 
\left( [a_{i-1}, \dots , a_1]+[a_{i+1}, \dots a_t] \right) \right). 
$$ 
Then, expanding into Taylor series according to 
$$
\frac{v}{n}-\frac{1}{n} \left( [a_{i-1}, \dots , a_1]+[a_{i+1}, \dots a_t] \right)
$$
function
$$
\frac{1}{q \left( a \right)^{2\beta}}, 
$$
when 
$$
|\frac{v}{n}-\frac{1}{n} \left( [a_{i-1}, \dots , a_1]+[a_{i+1}, \dots a_t] \right)|<1
$$
we get absolutely converging series
\begin{equation}
\label{2nd}
\begin{split}
\frac{1}{q (a)^{2\beta}}=
\frac{1}{n^{2\beta}\langle a_1, \dots , a_{j-1}\rangle^{2\beta}\langle a_{j+1}, \dots , a_t \rangle^{2\beta}}\cdot\\
\cdot\left( 1+
\sum_{k=1}^{\infty}\frac{2\beta \cdots \left( 2\beta+k-1 \right) }{k!}
 \left( \frac{v}{n}-\frac{1}{n}( [a_{i-1}, \dots , a_1]+[a_{i+1}, \dots a_t]) \right)^k \right)={}\\
=\frac{1}{n^{2\beta}\langle a_1, \dots , a_{j-1}\rangle^{2\beta}\langle a_{j+1}, \dots , a_t \rangle^{2\beta}}+{}\\
+\frac{1}{n^{2\beta}}\frac{1}{\langle a_1, \dots , a_{j-1}\rangle^{2\beta}\langle a_{j+1}, \dots , a_t \rangle^{2\beta}}
\sum_{k=1}^{\infty}\frac{1}{n^{k}}\gamma_k(2\beta)
 \left( v- ( [a_{i-1}, \dots , a_1]+[a_{i+1}, \dots a_t]) \right) ^k, 
\end{split}
\end{equation}
where
\begin{equation}
\label{gamma}
\gamma_k(\beta)=\frac{\beta \cdots (\beta+k-1)}{k!}. 
\end{equation}
After substituting (\ref{2nd}) to (\ref{1st}) regarding to lemma 6 with the given $ w $, replacing sum according to $ a_i $ with sum
 according to $ v $, we get 
\begin{equation}
\label{as_1} 
\begin{split}
\sum_{
\begin{array}{c}
a \in A_n, \\
\exists j: a_j > n-w
\end{array}
}
\frac{1}{q^{2\beta}}
= 
\sum_{v=1}^{w-1}
\sum_{u+s=v}
\sum_{
\begin{array}{c}
a \in A_n, \\
a_1+ \dots +a_{j-1}=u, \\
a_{j+1}+ \dots +a_t=s
\end{array}
}
\frac{1}{n^{2\beta}\langle a_1, \dots , a_{j-1}\rangle^{2\beta}\langle a_{j+1}, \dots , a_t \rangle^{2\beta}}
+ \\
%%$$  
%%$$  
+
\frac{1}{n^{2\beta}}
\sum_{v=1}^{w-1}
\sum_{u+s=v}
\sum_{
\begin{array}{c}
a \in A_n, \\
a_1+ \dots +a_{j-1}=u, \\
a_{j+1}+ \dots +a_t=s
\end{array}
}
\frac{1}{\langle a_1, \dots , a_{j-1}\rangle^{2\beta}\langle a_{j+1}, \dots , a_t \rangle^{2\beta}}
\sum_{k=1}^{\infty}\frac{1}{n^{k}}\gamma_k(2\beta)\cdot \\
%%  
%%  
%%$$  
%%$$  
\cdot \left( v- \left( [a_{i-1}, \dots , a_1]+[a_{i+1}, \dots a_t] \right) \right)^k.
\end{split} \end{equation}
Let's consider the main term in asymptotic formula we've got. 
According to lemma 7, 
$$ 
\frac{1}{n^{2\beta}}
\sum_{u+v \leq w}
\sum_{a_1+. . +a_{k}=u}\frac{1}{\langle a_1, \dots , a_k\rangle^{2\beta}}
\sum_{a_1+. . +a_{l}=v, a_l \ge 2}\frac{1}{\langle a_1, \dots , a_l\rangle^{2\beta}}=
$$ 
$$ 
=
\frac{C_0}{n^{2\beta}}
 +O \left( \frac{1}{n^{2\beta}w^{2 \left( \beta-1 \right)}} \right). 
$$ 
Now let's consider the error term:
$$ 
\frac{1}{n^{2\beta}}
\sum_{k=1}^{\infty}\frac{1}{n^{k}}
\cdot
\gamma_k(2\beta)
\sum_{v=1}^{w-1}
\sum_{
\begin{array}{c}
a \in A_n, \\
a_1+ \dots +a_{j-1}+a_{j+1}+ \dots +a_t=v
\end{array}
}
$$ 
$$ 
\frac{1}{\langle a_1, \dots , a_{j-1}\rangle^{2\beta}\langle a_{j+1}, \dots , a_t \rangle^{2\beta}}
 \left( v- \left( [a_{i-1}, \dots , a_1]+[a_{i+1}, \dots a_t] \right) \right)^k
+
O \left( \frac{1}{n^{2r \left( \beta-1 \right)}} \right)
=
$$  
$$
=
\frac{1}{n^{2\beta}}
\sum_{k=1}^{\infty}\frac{R_k}{n^{k}}
+
O \left( \frac{1}{n^{2r \left( \beta-1 \right)}} \right). 
$$ 
Coefficient $ R_{k} $ at $ k $th term is equal to
$$ 
R_{k}=\gamma_k(2\beta)\cdot
\sum_{v=1}^{w-1}
\sum_{
\begin{array}{c}
a \in A_n, \\
a_1+ \dots +a_{j-1}+a_{j+1}+ \dots +a_t=v
\end{array}
}
$$ 
$$ 
\frac{1}{\langle a_1, \dots , a_{j-1}\rangle^{2\beta}\langle a_{j+1}, \dots , a_t \rangle^{2\beta}}
 \left( v- \left( [a_{i-1}, \dots , a_1]+[a_{i+1}, \dots a_t] \right) \right)^k=
$$ 
$$ 
= \left( \sum_{v=1}^{\infty}-\sum_{v=w}^{\infty} \right)
\sum_{
\begin{array}{c}
a \in A_n, \\
a_1+ \dots +a_{j-1}+a_{j+1}+ \dots +a_t=v
\end{array}
}
\gamma_k(2\beta)\cdot
$$ 
$$ 
\cdot\frac{1}{\langle a_1, \dots , a_{j-1}\rangle^{2\beta}\langle a_{j+1}, \dots , a_t \rangle^{2\beta}}
 \left( v- \left( [a_{i-1}, \dots , a_1]+[a_{i+1}, \dots a_t] \right) \right)^k. 
$$ 

Let's investigate the convergence of series 
\begin{equation}
\label{series1}
\begin{split}
\sum_{v=1}^{\infty}
\sum_{
\begin{array}{c}
a \in A_n, \\
a_1+ \dots +a_{j-1}+a_{j+1}+ \dots +a_t=v
\end{array}
}
\frac{ \left( v- \left( [a_{i-1}, \dots , a_1]+[a_{i+1}, \dots a_t] \right) \right)^k}
{\langle a_1, \dots , a_{j-1}\rangle^{2\beta}\langle a_{j+1}, \dots , a_t \rangle^{2\beta}}.
\end{split}
\end{equation}
These series can be ameliorated with the series like 
$$ 
K_k\sum_{v=1}^{\infty}v^k
\sum_{
\begin{array}{c}
a \in A_n, \\
a_1+ \dots +a_{j-1}+a_{j+1}+ \dots +a_t=v
\end{array}
}
\frac{1}{\langle a_1, \dots , a_{j-1}\rangle^{2\beta}\langle a_{j+1}, \dots , a_t \rangle^{2\beta}}=
$$ 
$$ 
=K_k\sum_{v=1}^{\infty}v^k
\sum_{s+t=v}
\sum_{a_1+ \dots +a_j=s}
\frac{1}{\langle a_1, \dots , a_j\rangle^{2\beta}}
\sum_{a_1+ \dots +a_i=t}
\frac{1}{\langle a_1, \dots , a_i\rangle^{2\beta}}. 
$$ 
where $ K_k $ are some constants. 
According to lemma 9, 
$$ 
\sum_{a_1+ \dots +a_j=s}
\frac{1}{\langle a_1, \dots , a_j\rangle^{2\beta}}=O \left( \frac{1}{s^{2\beta}} \right), 
$$ 
$$ 
\sum_{a_1+ \dots +a_i=t}
\frac{1}{\langle a_1, \dots , a_i\rangle^{2\beta}}=O \left( \frac{1}{t^{2\beta}} \right), 
$$ 
hence, main term of series (\ref{series1}) can be estimated as $ O \left( \frac{1}{v^{2\beta-k-1}} \right) $. 
Thus, series converges when $ 2\beta-k-1 > 1 $, e. i. $ k<2\beta-2 $. \\
When $ k < 2\beta -2 $ constants $ {C'}_k $ can be defined as follows:
$$ 
{C'}_k=
\gamma_k(2\beta)
\sum_{v=1}^{\infty}
\sum_{a \in A_n, a_1+ \dots +a_{j-1}+a_{j+1}+ \dots +a_t=v}
\frac{\left( v- \left( [a_{i-1}, \dots , a_1]+[a_{i+1}, \dots a_t] \right) \right)^k }
{\langle a_1, \dots , a_{j-1}\rangle^{2\beta}\langle a_{j+1}, \dots , a_t \rangle^{2\beta}}.
$$ 
When $ k < 2\beta -2 $, let's estimate diversity 
$$ 
{C'}_k-\gamma_k(2\beta)
\sum_{v=1}^{w-1}
\sum_{
\begin{array}{c}
a \in A_n, \\
a_1+ \dots +a_{j-1}+a_{j+1}+ \dots +a_t=v
\end{array}
}
\frac{\left( v- \left( [a_{i-1}, \dots , a_1]+[a_{i+1}, \dots a_t] \right) \right)^k}
{\langle a_1, \dots , a_{j-1}\rangle^{2\beta}\langle a_{j+1}, \dots , a_t \rangle^{2\beta}}
=
$$
$$
=
\gamma_k(2\beta)
\sum_{v=w}^{\infty}
\sum_{
\begin{array}{c}
a \in A_n, \\
a_1+ \dots +a_{j-1}+a_{j+1}+ \dots +a_t=v
\end{array}
}
\frac{ \left( v- \left( [a_{i-1}, \dots , a_1]+[a_{i+1}, \dots a_t] \right) \right)^k}
{\langle a_1, \dots , a_{j-1}\rangle^{2\beta}\langle a_{j+1}, \dots , a_t \rangle^{2\beta}}
 \ll
$$ 
$$ 
\ll 
\sum_{v=w}^{\infty}v^k
\sum_{s+t=v}
\sum_{a_1+ \dots +a_j=s}
\frac{1}{\langle a_1, \dots , a_j\rangle^{2\beta}}
\sum_{a_1+ \dots +a_i=t}
\frac{1}{\langle a_1, \dots , a_i\rangle^{2\beta}} \ll
$$ 
$$ 
\ll
\int_{w}^{\infty}
\frac{d v}{v^{2\beta-k-1}}=O \left( \frac{1}{w^{2\beta-k-2}} \right). 
$$ 
Thus we obtain, that $ k $th term when $ k < 2\beta -2 $ is equal to
$$ 
\frac{\gamma_k(2\beta)}{n^{2\beta+k}}
\sum_{v=1}^{w-1}
\sum_{
\begin{array}{c}
a \in A_n, \\
a_1+ \dots +a_{j-1}+a_{j+1}+ \dots +a_t=v
\end{array}
}
\frac{\left( v- \left( [a_{i-1}, \dots , a_1]+[a_{i+1}, \dots a_t] \right) \right)^k}
{\langle a_1, \dots , a_{j-1}\rangle^{2\beta}\langle a_{j+1}, \dots , a_t \rangle^{2\beta}}
 =
$$ 
$$ 
=C'_k\frac{1}{n^{2\beta+k}}+O \left( \frac{1}{n^{2\beta+k}w^{2\beta-k-2}} \right)=
C'_k\frac{1}{n^{2\beta+k}}+O \left( \frac{1}{n^{2\beta}w^{2\beta-2}} \right). 
$$ 

Now let us consider the error term of the series in case $ k\ge 2\beta-2 $. 

When $ k>2\beta-2 $ 
$$ 
R_{k}=\gamma_k(2\beta)
\sum_{v=1}^{w-1}
\sum_{
\begin{array}{c}
a \in A_n, \\
a_1+ \dots +a_{j-1}+a_{j+1}+ \dots +a_t=v
\end{array}
}
\frac{\left( v- \left( [a_{i-1}, \dots , a_1]+[a_{i+1}, \dots a_t] \right) \right)^k}
{\langle a_1, \dots , a_{j-1}\rangle^{2\beta}\langle a_{j+1}, \dots , a_t \rangle^{2\beta}}
 \leq
$$
$$
\leq
\gamma_k(2\beta)
\sum_{v=1}^{w-1}v^k
\sum_{s+u=v}
\sum_{a_1+ \dots +a_j=s}
\frac{1}{\langle a_1, \dots , a_j\rangle^{2\beta}}
\sum_{a_1+ \dots +a_i=u}
\frac{1}{\langle a_1, \dots , a_i\rangle^{2\beta}} \leq
$$ 
$$ 
\leq
\gamma_k(2\beta)\cdot
16C_0^2
\int_{1}^{w-1}
\frac{d v}{v^{2\beta-k-1}}\leq
$$ 
$$ 
\leq
\gamma_k(2\beta)\cdot
16C_0^2
w^{k+2-2\beta}. 
$$ 
Then summing according to $ k > 2\beta-2 $, we get . 
$$ 
\sum_{k>2\beta-2}
\frac{R_{k}}{n^{2\beta+k}}
\leq
\sum_{k>2\beta-2}
\gamma_k(2\beta)\cdot
16C_0^2
\frac{w^{k+2-2\beta}}{n^{2\beta+k}}=
$$ 
$$ 
=
\frac{1}{n^{2\beta}w^{2\beta-2}}\cdot
16\left( C_0 \right)^2
\sum_{k>2\beta-2}
\gamma_k(2\beta)\cdot
\left(\frac{w}{n}\right)^k <
$$ 
$$ 
<
\frac{1}{n^{2\beta}w^{2\beta-2}}\cdot
16C_0^2
\sum_{k=1}^{\infty}
\gamma_k(2\beta)\cdot
\left(\frac{w}{n}\right)^k =
\frac{1}{n^{2\beta}w^{2\beta-2}}\cdot
16C_0^2
\left(\frac{1}{1-\frac{w}{n}}\right)^{2\beta}. 
$$ 
With the given $ w $ value of $ \frac{1}{1-\frac{w}{n}} $ doesn't emceed 2, 
 hence the sum can be estimated as $ O\left(\frac{1}{n^{2\beta}w^{2\beta-2}}\right) $. 
$ 2\beta $ can turn out to be integer.
In this case when $ k=2\beta-2 $ we get
$$ 
\gamma_k(2\beta)
\sum_{v=1}^{w-1}
\sum_{
\begin{array}{c}
a \in A_n, \\
a_1+ \dots +a_{j-1}+a_{j+1}+ \dots +a_t=v
\end{array}
}
\frac{1}{\langle a_1, \dots , a_{j-1}\rangle^{2\beta}\langle a_{j+1}, \dots , a_t \rangle^{2\beta}}
\times
$$ 
$$ 
\times
 \left( v- \left( [a_{i-1}, \dots , a_1]+[a_{i+1}, \dots a_t] \right) \right)^k=
O \left( \log w \right). 
$$ 

Thus, the error term in case when $ 2\beta $ is integer is
$$ 
O \left( \frac{\log w}{n^{4\beta-2}}+\frac{1}{w^{2\beta-2}n^{2\beta}} \right), 
$$ 
otherwise it is
$$ 
O \left( \frac{1}{w^{2\beta-2}n^{2\beta}} \right). 
$$ 
Thus, for $ \sum_{n, 1}^{ \left( 2 \right)} $ the following asymptotic holds: when $ 2\beta $ is integer 
$$ 
\Sigma_{n, 1}^{ \left( 2 \right)}
=
\frac{C_0}{n^{2\beta}}
+
O \left( \frac{1}{n^{2\beta}w^{2\beta-2}} \right)
+
$$ 
$$ 
+\sum_{1\le k<2\beta-2} \left( {C'}_k\frac{1}{n^{2\beta+k}}+O \left( \frac{1}{n^{2\beta}w^{2\beta-2}} \right) \right)
+
$$ 
$$ 
+
O \left( 
\frac{\log n}{n^{4\beta-2}} 
+
\frac{1}{n^{2\beta}w^{2 \left( \beta-1 \right)}}
+
\frac{1}{n^{2r \left( \beta-1 \right)}}
 \right)=
$$ 
$$ 
=
\frac{C_0}{n^{2\beta}}
+
\sum_{1\le k<2\beta-2}{C'}_k\frac{1}{n^{2\beta+k}}
+
O \left( 
\frac{\log n}{n^{4\beta-2}} 
+
\frac{1}{n^{2\beta}w^{2 \left( \beta-1 \right)}}
+
\frac{1}{n^{2r \left( \beta-1 \right)}}
 \right)
$$ 
when $ 2\beta $ is not integer , 
$$ 
\Sigma_{n, 1}^{ \left( 2 \right)}
=
\frac{C_0}{n^{2\beta}}
+
O \left( \frac{1}{n^{2\beta}w^{2\beta-2}} \right)
+
$$ 
$$ 
+\sum_{1\le k<2\beta-2} \left( {C'}_k\frac{1}{n^{2\beta+k}}+O \left( \frac{1}{n^{2\beta}w^{2\beta-2}} \right) \right)
+
O \left( 
\frac{1}{n^{2\beta}w^{2 \left( \beta-1 \right)}}
+
\frac{1}{n^{2r \left( \beta-1 \right)}}
 \right)=
$$ 
$$ 
=
\frac{C_0}{n^{2\beta}}
+
\sum_{1\le k<2\beta-2}{C'}_k\frac{1}{n^{2\beta+k}}
+
O \left( 
\frac{1}{n^{2\beta}w^{2 \left( \beta-1 \right)}}
+
\frac{1}{n^{2r \left( \beta-1 \right)}}
 \right).
$$ 

{\bf Lemma is proved. }
\newline

{\bf Theorem is proved. }

\medskip
\medskip
{\large \bf\S4. Proof of the main result. }\\

\medskip
Let's remind, that $ r $ and $ w $ are integer parameters, satisfying conditions
 $ r \geq 1 $ and $ 1 \leq w \leq n $. \\
We'll use the partition of $ A_n $, 
defined in the proof of theorem 3. 
Then, we divide $ A_{n, 1}^{ \left( 2 \right)} $ into 2 sets: $ A_{n, 1}^{ \left( 3 \right)} $, where the greatest partial quotient is the last one, 
and $ A_{n, 2}^{ \left( 3 \right)} $, where it's not the last one: 
$$ 
A_{n, 1}^{ \left( 3 \right)} = \{ a \in A_{n, 1}^{ \left( 2 \right)}
 \mid a_t > \max \{ a_1, \ldots, a_{t-1} \} \}
$$ 
and
$$ 
A_{n, 2}^{ \left( 3 \right)} = A_{n, 1}^{ \left( 2 \right)} \backslash A_{n, 1}^{ \left( 3 \right)} = 
\{ a = \left( a_1, \ldots, a_t \right) \in A_{n, 1}^{ \left( 2 \right)}
 \mid a_t \leq \max \{ a_1, \ldots, a_{t-1} \} \}
.$$ 
For subset $ A_{n, i}^{ \left( j \right)} $ of $ A $ let's define
$$ 
\Sigma_{n, i}^{ \left( j \right)} = \sum_{a \in A_{n, i}^{ \left( j \right)}}
 \frac{1} { \left( qq_{-} \right)^{\beta}}+
 \frac{1} { \left( qq_{+} \right)^{\beta}}, 
$$ 
where $ a= \left( a_1, \dots , a_t \right), q=q \left( a \right)=\langle a_1, \dots , a_t \rangle $ ; $ q_{-}=q_{-} \left( a \right) $ and $ q_{+}=q_{+} \left( a \right) $ 
are defined in lemma 1. 
\newline
Then we divide $ \sum_{n, 1}^{ \left( 3 \right)} $ into $ \sum_{n, 1}^{ \left( 3 \right)+} $ and $ \sum_{n, 1}^{ \left( 3 \right)-} $ with

$$ 
\Sigma_{n, 1}^{ \left( 3 \right)+} = \sum_{a \in A_{n, 1}^{ \left( 3 \right)}}
 \frac{1} { \left( qq_{+} \right)^{\beta}}, 
\Sigma_{n, 1}^{ \left( 3 \right)-} = \sum_{a \in A_{n, 1}^{ \left( 3 \right)}}
 \frac{1} { \left( qq_{-} \right)^{\beta}}. 
$$

Thus, 
\begin{equation}
\label{whole_sum}
\sigma \left( F_n \right) = 
\Sigma_{n, 2}^{ \left( 1 \right)} +
\Sigma_{n, 2}^{ \left( 2 \right)} +
\Sigma_{n, 2}^{ \left( 3 \right)} +
\Sigma_{n, 1}^{ \left( 3 \right)+} + \Sigma_{n, 1}^{ \left( 3 \right)-}
. 
\end{equation}
\\
Let's estimate these sums separately. \\
According to [5], 
for $ \sum_{n, 2}^{ \left( 1 \right)} $ 
the following estimate holds:
\newline
{\bf Lemma 11. }
{\it
\begin{equation}
\label{big_q}
\Sigma_{n, 2}^{ \left( 1 \right)}
\le
\frac{1}{n^{ \left( \beta-1 \right) \left( 2r-1 \right)}}. 
\end{equation}
}
\newline

{\bf Lemma 12. }\\
{\it
When $ n \rightarrow \infty $ we have
\begin{equation}
\label{sigma_n22}
\Sigma_{n, 2}^{ \left( 2 \right)}
\ll
\frac{n^2\log^{3\beta}n}{w^{3\beta}}. 
\end{equation}
}
Lemma 12 is an analogue of lemma 5. \\
{\bf Proof. }
\\
According to lemma 4 for every $ a \in A_{n, 1}^{ \left( 1 \right)} $ holds $ t \le K r \log n $. 
As $ n = a_1 + \ldots + a_t \leq t \max a_j $, then
 $ \max a_j \geq \frac{n}{K r \log n} $. 
\\
Let $ a \in A_{n, 2}^{ \left( 2 \right)} $ and $ j $ be such, that
 $ a_j = \max \{ a_1, \ldots, a_t\}. $ 
As $ a_j \leq n - w $, then for the sum of other $ a_j $ we have 
 $ \sum\limits_{i \not= j} a_j \ge w $, and, similarly to the above, we have
 $ \max_{i \not= j} a_j \geq \frac{w}{K r \log n} $. 
\newline
Thus, there's at least one index $ j \leq t-1 $ such, 
that $ a_j \geq \frac{w}{K r \log n} $. 
Hence, 
$$ 
\Sigma_{n, 2}^{ \left( 2 \right)} \le
\sum_{
\begin{array}{c}
a_1+\ldots+a_t = n, \\
\langle a_1, \dots , a_t \rangle \leq n^r, \\
\exists k, l, k\neq l: a_k, a_l \geq \frac{w}{K\log n}
\end{array}
}
\left(
\frac{1}{ \left( qq_{-} \right)^{\beta}}+
\frac{1}{ \left( qq_{+} \right)^{\beta}} 
\right)
\leq
$$ 
$$ 
\leq
\sum_{
\begin{array}{c}
a_1+\ldots+a_t = n, \\
\langle a_1, \dots , a_t \rangle \leq n^r, \\
\exists k \leq \left( t-1 \right) : a_k, a_t \geq \frac{w}{K\log n}
\end{array}
}
\left(
\frac{1}{ \left( qq_{-} \right)^{\beta}}+
\frac{1}{ \left( qq_{+} \right)^{\beta}}
\right)
+
$$
$$
+
\sum_{
\begin{array}{c}
a_1+\ldots+a_t = n, \\
\langle a_1, \dots , a_t \rangle \leq n^r, \\
\exists k,l \leq \left( t-1 \right) : a_k, a_l \geq \frac{w}{K\log n}
\end{array}
}
\left(
\frac{1}{ \left( qq_{-} \right)^{\beta}}+
\frac{1}{ \left( qq_{+} \right)^{\beta}}
\right)
\leq
$$ 
$$ 
\leq
2\sum_{
\begin{array}{c}
a_1+\ldots+a_t = n, \\
\langle a_1, \dots , a_t \rangle \leq n^r, \\
\exists k \leq \left( t-1 \right) : a_k, a_t \geq \frac{w}{K\log n}
\end{array}
}
\frac{1}{ \left( qq_{-} \right)^{\beta}}+ 
$$
$$
+2\sum_{
\begin{array}{c}
a_1+\ldots+a_t = n, \\
\langle a_1, \dots , a_t \rangle \leq n^r, \\
\exists k,l \leq \left( t-1 \right) : a_k, a_l \geq \frac{w}{K\log n}
\end{array}
}
\frac{1}{ \left( qq_{-} \right)^{\beta}}. 
$$ 
Using (\ref{knut}), we obtain the following estimation for the continuant:
$$ 
q \left( a \right) \geq a_k a_t\langle a_1, . . , a_{k-1}\rangle
\langle a_{k+1}, . . , a_{t-1}\rangle 
$$ 
and
$$
q \left( a \right) \geq a_k a_l\langle a_1, . . , a_{k-1}\rangle
\langle a_{k+1}, . . , a_{l-1}\rangle 
\langle a_{l+1}, . . , a_{t}\rangle 
$$
Thus, splitting sum $\sum_{n,2}^(2)$ into two parts (one part corresponds to items with big last patial quotient $a_t$, and another part corresponds to items which have big partial quotints $a_k,a_l$,such that neither of them is the last one), we have
%%$$ 
\begin{equation}
\label{sum_2parts}
\begin{split}
\Sigma_{n, 2}^{ \left( 2 \right)}
\ll
\sum_{
\begin{array}{c}
a_k+a_t \le n, \\
a_k, a_t \geq \frac{w}{K\log n}
\end{array}
}
\frac{1}{a_k^{2\beta}a_t^{\beta}}
\sum_{u+v=n-a_k-a_t}  \\
%%$$
%%$$ 
\sum_{a_1+ \dots +a_{k-1}=u}\frac{1}{\langle a_1, \dots , a_{k-1} \rangle^{2\beta}}
\sum_{a_{k+1}+ \dots +a_{t-1}=v}\frac{1}{\langle a_{k+1}, \dots , a_{t-1} \rangle^{2\beta}}\\
%\ll
%%$$
%%$$  
+
\sum_{
\begin{array}{c}
a_k+a_l \le n, \\
a_k, a_l \geq \frac{w}{K\log n}
\end{array}
}
\frac{1}{a_k^{2\beta}a_l^{2\beta}}
\sum_{
u+v+s=n-a_k-a_l
}
%%$$  
%%$$  
\sum_{a_1+ \dots +a_{k-1}=u}\frac{1}{\langle a_1, \dots , a_{k-1} \rangle^{2\beta}}\\
\sum_{a_{k+1}+ \dots +a_{l-1}=v}\frac{1}{\langle a_{k+1}, \dots , a_{l-1} \rangle^{2\beta}} 
\sum_{a_{k+1}+ \dots +a_{l-1}=s}\frac{1}{\langle a_{l+1}, \dots , a_{t-1} \rangle^{\beta}\langle a_{l+1}, \dots , a_t \rangle^{\beta}}. 
%\ll
\end{split}
\end{equation}  
Let's consider the inner sum in the first item in (\ref{sum_2parts}) . 
$$ 
\sum_{u+v=n-a_k-a_t}
\sum_{a_1+ \dots +a_{k-1}=u}\frac{1}{\langle a_1, \dots , a_{k-1} \rangle^{2\beta}}
\sum_{a_{k+1}+ \dots +a_{t-1}=v}\frac{1}{\langle a_{k+1}, \dots , a_{t-1} \rangle^{2\beta}}=
$$ 
$$ 
\sum_{u+v=n-a_k-a_t}
\sum_{x_1+ \dots +x_r=u}\frac{1}{\langle x_1, \dots , x_r \rangle^{2\beta}}
\sum_{y_1+ \dots +y_h=v}\frac{1}{\langle y_1, \dots , y_h \rangle^{2\beta}}\leq
$$ 
$$ 
\leq
\sum_{u+v\leq n}
\sum_{x_1+ \dots +x_r=u}\frac{1}{\langle x_1, \dots , x_r \rangle^{2\beta}}
\sum_{y_1+ \dots +y_h=v}\frac{1}{\langle y_1, \dots , y_h \rangle^{2\beta}}\leq
$$ 
$$ 
\leq
\sum_{u+v\leq \infty}
\sum_{x_1+ \dots +x_r=u}\frac{1}{\langle x_1, \dots , x_r \rangle^{2\beta}}
\sum_{y_1+ \dots +y_h=v}\frac{1}{\langle y_1, \dots , y_h \rangle^{2\beta}}\leq
$$ 
$$ 
\leq
\left(1+2\sum_{
\begin{array}{c}
x_1+ \dots +x_r \leq \infty, \\
x_r \geq 2
\end{array}
}
\frac{1}{\langle x_1, \dots , x_r \rangle^{2\beta}}\right)^2
\leq
$$ 
$$ 
\leq 
\left(1+
2\frac{\zeta(2\beta-1)}{\zeta{2\beta}}
\right)^2. 
$$ 
Then for the outer sum for the fist item in (\ref{sum_2parts}) holds
$$ 
\sum_{
\begin{array}{c}
a_k+a_t \le n, \\
a_k, a_t \geq \frac{w}{K\log n}
\end{array}
}
\frac{1}{ \left( a_k \right)^{2\beta} \left( a_t \right)^{\beta}}
\ll
\frac{\log^{3\beta}n}{w^{3\beta}}
 \left( n-2\frac{w}{K \log n} \right)^2. 
$$ 
Here $ \left( n-2\frac{w}{K \log n} \right)^2 $ is the number of elements in sum, 
 $ \frac{\log^{3\beta}n}{w^{3\beta}} $ -- upper bound for the value under summation. 
%%2nd sum 
Now let's consider the inner sum in the second item in (\ref{sum_2parts}). 
$$  
\sum_{u+v+s=n-a_k-a_l}
\sum_{a_1+ \dots +a_{k-1}=u}\frac{1}{\langle a_1, \dots , a_{k-1} \rangle^{2\beta}}
\times
$$  
$$  
\times
\sum_{a_{k+1}+ \dots +a_{t-l}=v}\frac{1}{\langle a_{k+1}, \dots , a_{t-1} \rangle^{2\beta}}
\times
$$  
$$  
\times
\sum_{a_{k+1}+ \dots +a_{l-1}=s}\frac{1}{\langle a_{l+1}, \dots , a_{t-1} \rangle^{\beta}\langle a_{l+1}, \dots , a_t \rangle^{\beta}}
=
$$  
$$  
=
\sum_{u+v+s=n-a_k-a_l}
\sum_{x_1+ \dots +x_j=u}\frac{1}{\langle x_1, \dots , x_j \rangle^{2\beta}}
\times
$$
$$
\times
\sum_{y_1+ \dots +y_h=v}\frac{1}{\langle y_1, \dots , y_h \rangle^{2\beta}}
\sum_{z_1+ \dots +z_p=s}\frac{1}{\langle z_1, \dots , z_{p-1} \rangle^{\beta}\langle z_{1}, \dots , z_p \rangle^{\beta}}
\leq
$$  
$$  
\leq
\sum_{u+v+s\leq n}
\sum_{x_1+ \dots +x_j=u}\frac{1}{\langle x_1, \dots , x_j \rangle^{2\beta}}
\sum_{y_1+ \dots +y_h=v}\frac{1}{\langle y_1, \dots , y_h \rangle^{2\beta}}
\times
$$  
$$  
\times
\sum_{z_1+ \dots +z_p=s}\frac{1}{\langle z_1, \dots , z_{p-1} \rangle^{\beta}\langle z_{1}, \dots , z_p \rangle^{\beta}}
\leq
$$  
$$  
\leq
\sum_{u+v+s\leq \infty}
\sum_{x_1+ \dots +x_j=u}\frac{1}{\langle x_1, \dots , x_j \rangle^{2\beta}}
\sum_{y_1+ \dots +y_h=v}\frac{1}{\langle y_1, \dots , y_h \rangle^{2\beta}}
\times
$$  
$$  
\times
\sum_{z_1+ \dots +z_p=s}\frac{1}{\langle z_1, \dots , z_{p-1} \rangle^{\beta}\langle z_{1}, \dots , z_p \rangle^{\beta}}
\leq
$$  
$$  
\leq
\left(1+2\sum_{
\begin{array}{c}
x_1+ \dots +x_j \leq \infty, \\
x_j \geq 2
\end{array}
}
\frac{1}{\langle x_1, \dots , x_j \rangle^{2\beta}}\right)^2 
\times
$$
$$
\times
\left(1+\sum_{
\begin{array}{c}
z_1+ \dots +z_p \leq \infty, \\
z_p \geq 2
\end{array}
}
\frac{1}{\langle z_1, \dots , z_{p-1} \rangle^{\beta}\langle z_1, \dots , z_{p} \rangle^{\beta}}\right)
\leq
$$  
$$  
\leq 
\left(1+
2\frac{\zeta(2\beta-1)}{\zeta(2\beta)}
\right)^2
\left(1+
\frac{\zeta(2\beta-1)}{\zeta(2\beta)}
\right)
. 
$$  
Then for the outer sum in the second item in (\ref{sum_2parts}) the following estimation holds
$$  
\sum_{
\begin{array}{c}
a_k+a_t \le n, \\
a_k, a_l \geq \frac{w}{K\log n}
\end{array}
}
\frac{1}{ \left( a_k \right)^{2\beta} \left( a_l \right)^{2\beta}}
\ll
\frac{\log^{4\beta}n}{w^{4\beta}}
 \left( n-2\frac{w}{K \log n} \right)^2. 
$$  
Here $ \left( n-2\frac{w}{K \log n} \right)^2 $ is the number of elements in sum, 
 $ \frac{\log^{4\beta}n}{w^{4\beta}} $ -upper bound for the value under summation. \\
Thus, for $ \sum_{n, 2}^{(2)} $ we get
$$ 
\Sigma_{n, 2}^{(2)}
=
O \left( \frac{n^2\log^{3\beta}n}{w^{3\beta}}
 \right). 
$$ 
{\bf Lemma is proved. }
\\

{\bf Lemma 13. }\\
{\it
When $ w \leq \frac{n}{2}-2 $ the following asymptotic holds for the sum $ \sum_{n, 2}^{(3)} $:\\
in the case $ \beta $ is integer, 
\begin{equation}
\label{sigma_n23_int}
\Sigma_{n, 2}^{ \left( 3 \right)}
=
\sum_{0\le k<\beta-2}
B_k
\frac{1}{n^{2\beta+k}}
+O\left(
\frac{1}{n^{(\beta-1)(2r-1)}}
+
\frac{1}{n^{2\beta}w^{\beta-2}}
+
\frac{\log w}{n^{3\beta-2}}
\right), 
\end{equation}
otherwise, 
\begin{equation}
\label{sigma_n23}
\Sigma_{n, 2}^{ \left( 3 \right)}
=
\sum_{0\le k<\beta-2}
B_k
\frac{1}{n^{2\beta+k}}
+O\left(
\frac{1}{n^{(\beta-1)(2r-1)}}
+
\frac{1}{n^{2\beta}w^{\beta-2}}
\right), 
\end{equation}
where $ B_k $ are some constants. 
}
\\
{\bf Proof. }
$$ 
\Sigma_{n, 2}^{ \left( 3 \right)}
=
\sum_{
\begin{array}{c}
a \in A_n, \\
q \left( a \right) < n^r, \\
a_i > n-w, j \ne t
\end{array}
}
\left(
\frac{1}{ \left( qq_{-} \right)^{\beta}}
+
\frac{1}{ \left( qq_{+} \right)^{\beta}}
\right)=
$$ 
$$ 
=
\sum_{
\begin{array}{c}
a \in A_n, \\
a_i > n-w, j \ne t
\end{array}
}
\left(
\frac{1}{ \left( qq_{-} \right)^{\beta}}
+
\frac{1}{ \left( qq_{+} \right)^{\beta}}
\right)
-
\sum_{
\begin{array}{c}
a \in A_n, \\
q \left( a \right) \geq n^r, \\
a_i > n-w, j \ne t
\end{array}
}
\left(
\frac{1}{ \left( qq_{-} \right)^{\beta}}
+
\frac{1}{ \left( qq_{+} \right)^{\beta}}
\right). 
$$ 

The second sum can be estimated according to lemma 11 as 
$$ 
\sum_{\begin{array}{c}
a \in A_n, \\
q \left( a \right) \geq n^r, \\
a_i > n-w, j \ne t
\end{array}}
\left(
\frac{1}{ \left( qq_{-} \right)^{\beta}}
+
\frac{1}{ \left( qq_{+} \right)^{\beta}}
\right)
=
O \left( 
\frac{1}{n^{ \left( \beta-1 \right) \left( 2r-1 \right)}}
 \right). 
$$ 
Let us consider the first sum. 
Let $ a_i=n-v, v=1, \dots , (w-1) $. 
Using (\ref{knut}) we get
$$ 
\langle a_1, \dots , a_i, \dots , a_t\rangle=\langle a_1, \dots , a_{j-1}\rangle
\langle a_{j+1}, \dots , a_t\rangle
\left( a_i+[a_{j-1}, \dots , a_1]+[a_{j+1}, \dots , a_t]\right)=
$$ 
$$ 
=
n\langle a_1, \dots , a_{j-1}\rangle
\langle a_{j+1}, \dots , a_t\rangle
\left( 1-\frac{v}{n}+\frac{1}{n}\left([a_{j-1}, \dots , a_1]+[a_{j+1}, \dots , a_t]
\right)\right), 
$$ 
$$ 
\langle a_1, \dots , a_i, \dots , a_{t-1}\rangle
=
n\langle a_1, \dots , a_{i-1}\rangle
\langle a_{i+1}, \dots , a_{t-1}\rangle
\left( 1-\frac{v}{n}+\frac{1}{n}\left([a_{i-1}, \dots , a_1]+[a_{i+1}, \dots , a_{t-1}]
\right)\right), 
$$ 
$$ 
\langle a_1, \dots , a_i, \dots , a_t-1\rangle
=
n\langle a_1, \dots , a_{i-1}\rangle
\langle a_{i+1}, \dots , a_t-1\rangle
\cdot
$$
$$
\cdot
\left( 1-\frac{v}{n}+\frac{1}{n}\left([a_{i-1}, \dots , a_1]+[a_{i+1}, \dots , a_t-1]
\right)\right). 
$$ 
Then for $ \frac{1}{(qq_{-})^{\beta}} $ and $ \frac{1}{(qq_{+})^{\beta}} $ we can obtain the following formulas:
\begin{equation}
\label{ff1}
\begin{split}
\frac{1}{(qq_{-})^{\beta}}=
\frac{1}{n^{2\beta}}
\frac{1}{\langle a_1, \dots , a_{i-1}\rangle^{2\beta}}
\frac{1}{
\left(
\langle a_{i+1}, \dots , a_t\rangle
\langle a_{i+1}, \dots , a_{t-1}\rangle
\right)^{\beta}
}\cdot\\
\cdot
\frac{1}{
\left( 1-\frac{v}{n}+\frac{1}{n}\left([a_{i-1}, \dots , a_1]+[a_{i+1}, \dots , a_t]
\right)\right)^{\beta}
\left( 1-\frac{v}{n}+\frac{1}{n}\left([a_{i-1}, \dots , a_1]+[a_{i+1}, \dots , a_{t-1}]
\right)\right)^{\beta}
}, 
\end{split}
\end{equation}
\begin{equation}
\label{ff2}
\begin{split}
\frac{1}{(qq_{+})^{\beta}}=
\frac{1}{n^{2\beta}}
\frac{1}{\langle a_1, \dots , a_{i-1}\rangle^{2\beta}}
\frac{1}{
\left(
\langle a_{i+1}, \dots , a_t\rangle
\langle a_{i+1}, \dots , a_t-1\rangle
\right)^{\beta}
}
\cdot\\
\cdot
\frac{1}{
\left( 1-\frac{v}{n}+\frac{1}{n}\left([a_{i-1}, \dots , a_1]+[a_{i+1}, \dots , a_t]
\right)\right)^{\beta}
\left( 1-\frac{v}{n}+\frac{1}{n}\left([a_{i-1}, \dots , a_1]+[a_{i+1}, \dots , a_t-1]
\right)\right)^{\beta}
}. 
\end{split}
\end{equation}

Let's expand
$$ 
\frac{1}
{
\left(
1-\frac{v}{n}+\frac{1}{n}\left([a_{i-1}, \dots , a_1]+[a_{i+1}, \dots , a_t]\right)
\right)^{\beta}
}
$$ 
and 
$$ 
\frac{1}
{\left(1-\frac{v}{n}+\frac{1}{n}\left([a_{i-1}, \dots , a_1]+[a_{i+1}, \dots , a_{t-1}]\right)\right)^{\beta}}
$$ 
into Taylor series according to parameters 
$$ 
\frac{v}{n}-\frac{1}{n}\left([a_{i-1}, \dots , a_1]+[a_{i+1}, \dots , a_t]\right)
$$ 
and
$$ 
\frac{v}{n}-\frac{1}{n}\left([a_{i-1}, \dots , a_1]+[a_{i+1}, \dots , a_{t-1}]\right)
$$ 
correspondingly. Thus we obtain 
\begin{equation}
\label{fff1}
\begin{split}
\frac{1}
{\left(1-\frac{v}{n}+\frac{1}{n}\left([a_{i-1}, \dots , a_1]+[a_{i+1}, \dots , a_t]\right)\right)^{\beta}}
=\\
=
1
+
\sum_{k=1}^{\infty}\frac{1}{n^k}
\gamma_k(\beta)\cdot
\left(v-\left([a_{i-1}, \dots , a_1]+[a_{i+1}, \dots , a_t]\right)\right)^k
\end{split}
\end{equation}
when $ |\frac{v}{n}-\frac{1}{n}\left([a_{i-1}, \dots , a_1]+[a_{i+1}, \dots , a_t]\right)| <1 $, 
\begin{equation}
\label{fff2}
\begin{split}
\frac{1}
{\left(1-\frac{v}{n}+\frac{1}{n}\left([a_{i-1}, \dots , a_1]+[a_{i+1}, \dots , a_{t-1}]\right)\right)^{\beta}}=\\
=
1
+
\sum_{k=1}^{\infty}\frac{1}{n^k}
\gamma_k(\beta)
\cdot\left(v-\left([a_{i-1}, \dots , a_1]+[a_{i+1}, \dots , a_{t-1}]\right)\right)^k. 
\end{split}
\end{equation}
when $ |\frac{v}{n}-\frac{1}{n}\left([a_{i-1}, \dots , a_1]+[a_{i+1}, \dots , a_{t-1}]\right)| <1 $, 
where $ \gamma_k(\beta) $ are defined in (\ref{gamma}). 
When $ v \le n-1 $ series converge absolutely. 

Substituting (\ref{fff1}) and (\ref{fff2}) into (\ref{ff1}), we get with the given $ v $ 
\begin{equation}
\label{nonumber}
\begin{split}
%%$$  
\frac{1}{\left(qq_{-}\right)^{\beta}}=
\frac{1}{n^{2\beta}\langle a_1, \dots , a_{i-1}\rangle^{2\beta}\langle a_{i+1}, \dots , a_t\rangle^{\beta}
\langle a_{i+1}, \dots , a_{t-1}\rangle^{\beta}}\cdot \\
%%$$  
%%$$  
\cdot
\left( 1+ \right .
\sum_{k=1}^{\infty}\frac{1}{n^k}
%\cdot
\sum_{l+m=k}
\gamma_l(\beta)
\left(v-\left([a_{i-1}, \dots , a_1]+[a_{i+1}, \dots , a_{t-1}]\right)\right)^l \cdot \\
\cdot
\gamma_m(\beta)
\left(v-\left([a_{i-1}, \dots , a_1]+[a_{i+1}, \dots , a_t]\right)\right)^m
\left . \right). 
%%$$  
\end{split}
\end{equation}
Substituting the obtained result for $ \frac{1}{\left(qq_{-}\right)^{\beta}} $ into
$$ 
\sum_{v=1}^{w-1}
\sum_{a \in A_n, a_1+ \dots +a_{i-1}+a_{i+1}+ \dots +a_t=v}
\frac{1}{\left(qq_{-}\right)^{\beta}}=
$$ 
$$ 
=\sum_{v=1}^{w-1}
\sum_{a \in A_n, a_1+ \dots +a_{i-1}+a_{i+1}+ \dots +a_t=v}
\frac{1}{n^{2\beta}\langle a_1, \dots , a_{i-1}\rangle^{2\beta}\langle a_{i+1}, \dots , a_t\rangle^{\beta}
\langle a_{i+1}, \dots , a_{t-1}\rangle^{\beta}}\cdot
$$ 
$$ 
\cdot\left( \right .
1+\sum_{k=1}^{\infty}\frac{1}{n^k}
\sum_{l+m=k}
\gamma_l(\beta)
\left(v-\left([a_{i-1}, \dots , a_1]+[a_{i+1}, \dots , a_{t-1}]\right)\right)^l
\cdot
$$
$$
\cdot
\gamma_m(\beta)
\left(v-\left([a_{i-1}, \dots , a_1]+[a_{i+1}, \dots , a_t]\right)\right)^m
\left . \right)=
$$ 
$$ 
=\frac{1}{n^{2\beta}}
\sum_{v=1}^{w-1}
\sum_{a \in A_n, a_1+ \dots +a_{i-1}+a_{i+1}+ \dots +a_t=v}
\frac{1}{\langle a_1, \dots , a_{i-1}\rangle^{2\beta}\langle a_{i+1}, \dots , a_t\rangle^{\beta}
\langle a_{i+1}, \dots , a_{t-1}\rangle^{\beta}}
+
$$ 
$$ 
\sum_{k=1}^{\infty}
\frac{1}{n^{2\beta+k}}
\sum_{v=1}^{w-1}
\sum_{a \in A_n, a_1+ \dots +a_{i-1}+a_{i+1}+ \dots +a_t=v}
\frac{1}{\langle a_1, \dots , a_{i-1}\rangle^{2\beta}\langle a_{i+1}, \dots , a_t\rangle^{\beta}
\langle a_{i+1}, \dots , a_{t-1}\rangle^{\beta}}
\cdot
$$ 
$$ 
\cdot\sum_{l+m=k}
\gamma_l(\beta)
\left(v-\left([a_{i-1}, \dots , a_1]+[a_{i+1}, \dots , a_{t-1}]\right)\right)^l
%\cdot
%$$  
%$$  
\cdot
\gamma_m(\beta)
\left(v-\left([a_{i-1}, \dots , a_1]+[a_{i+1}, \dots , a_t]\right)\right)^m. 
$$ 
Let's investigate sum at $ \frac{1}{n^{2\beta+k}} $. 
$$ 
R_{k}^{-}=\sum_{v=1}^{w-1}
\sum_{a \in A_n, a_1+ \dots +a_{i-1}+a_{i+1}+ \dots +a_t=v}
\frac{1}{\langle a_1, \dots , a_{i-1}\rangle^{2\beta}\langle a_{i+1}, \dots , a_t\rangle^{\beta}
\langle a_{i+1}, \dots , a_{t-1}\rangle^{\beta}}\cdot
$$ 
$$ 
\cdot
\sum_{l+m=k}
\gamma_l(\beta)
\left(v-\left([a_{i-1}, \dots , a_1]+[a_{i+1}, \dots , a_{t-1}]\right)\right)^l
\cdot
\gamma_m(\beta)
\left(v- \left( [a_{i-1}, \dots , a_1]+[a_{i+1}, \dots , a_t]\right)\right)^m=
$$ 
$$ 
=
\left(\sum_{v=1}^{\infty}-\sum_{v=w}^{\infty}\right)
\sum_{a \in A_n, a_1+ \dots +a_{i-1}+a_{i+1}+ \dots +a_t=v}
\frac{1}{\langle a_1, \dots , a_{i-1}\rangle^{2\beta}\langle a_{i+1}, \dots , a_t\rangle^{\beta}
\langle a_{i+1}, \dots , a_{t-1}\rangle^{\beta}}\cdot
$$ 
$$ 
\cdot\sum_{l+m=k}
\gamma_l(\beta)
\left(v- \left( [a_{i-1}, \dots , a_1]+[a_{i+1}, \dots , a_{t-1}]\right)\right)^l
\cdot
\gamma_m(\beta)
\left(v- \left( [a_{i-1}, \dots , a_1]+[a_{i+1}, \dots , a_t]\right)\right)^m. 
$$ 
Investigate convergence of the first series. 
$$ 
\sum_{v=1}^{\infty}
\sum_{
\begin{array}{c}
a \in A_n, \\
a_1+ \dots +a_{i-1}+a_{i+1}+ \dots +a_t=v
\end{array}}
\frac{1}{\langle a_1, \dots , a_{i-1}\rangle^{2\beta}\langle a_{i+1}, \dots , a_t\rangle^{\beta}
\langle a_{i+1}, \dots , a_{t-1}\rangle^{\beta}}
$$ 
$$ 
\sum_{l+m=k}
 \gamma_l(\beta)
 \left( v- \left( [a_{i-1}, \dots , a_1]+[a_{i+1}, \dots , a_{t-1}]\right)\right)^l
\cdot
 \gamma_m(\beta)
 \left( v- \left( [a_{i-1}, \dots , a_1]+[a_{i+1}, \dots , a_t]\right)\right)^m \ll
$$ 
$$ 
\ll
\sum_{v=1}^{\infty}v^k
\sum_{\begin{array}{c}
a \in A_n, \\
a_1+ \dots +a_{i-1}+a_{i+1}+ \dots +a_t=v
\end{array}}
\frac{1}{\langle a_1, \dots , a_{i-1}\rangle^{2\beta}\langle a_{i+1}, \dots , a_t\rangle^{\beta}
\langle a_{i+1}, \dots , a_{t-1}\rangle^{\beta}} \ll
$$ 
$$ 
\ll\sum_{v=1}^{\infty}v^k
\sum_{\nu+\eta=v}
\sum_{a_1+ \dots +a_{i-1}=\nu}
\frac{1}{\langle a_1, \dots , a_{i-1}\rangle^{2\beta}}
\sum_{a_{i+1}+ \dots +a_t=\eta}
\frac{1}
{\langle a_{i+1}, \dots , a_t\rangle^{\beta}
\langle a_{i+1}, \dots , a_{t-1}\rangle^{\beta}}\ll
$$ 
$$ 
\ll
\sum_{v=1}^{\infty}v^k
\sum_{\nu+\eta=v}
\sum_{a_1+ \dots +a_t=\nu}
\frac{1}{q \left( a\right)^{2\beta}}
\sum_{a_1+ \dots +a_t'=\eta}
\frac{1}
{ \left( q \left( a\right)q_{-} \left( a\right)\right)^{\beta}}. 
$$ 
It follows from lemma 9, that 
$$ 
\sum_{a_1+ \dots +a_t=\nu}
\frac{1}{q \left( a\right)^{2\beta}}=O \left( \frac{1}{\nu^{2\beta}}\right), 
$$ 
it follows from lemma in [3] that
$$ 
\sum_{a_1+ \dots +a_t'=\eta}
\frac{1}
{ \left( q \left( a\right)q_{-} \left( a\right)\right)^{\beta}}=O \left( \frac{1}{\eta^{\beta}}\right). 
$$ 
Thus, common term of the series can be estimated as $ O \left( \frac{1}{v^{\beta-k-1}}\right) $, 
so the series converges when $ \beta-k-1 >1 $, i.e. when $ k<\beta-2 $. \\
With this $ k $ let us estimate the error term of the series:
$$ 
\sum_{v=w}^{\infty}
\sum_{
\begin{array}{c}
a \in A_n, \\
a_1+ \dots +a_{i-1}+a_{i+1}+ \dots +a_t=v
\end{array}
}
\frac{1}{\langle a_1, \dots , a_{i-1}\rangle^{2\beta}\langle a_{i+1}, \dots , a_t\rangle^{\beta}
\langle a_{i+1}, \dots , a_{t-1}\rangle^{\beta}}
$$ 
$$ 
\sum_{l+m=k}
 \gamma_l(\beta)
 \left( v- \left( [a_{i-1}, \dots , a_1]+[a_{i+1}, \dots , a_{t-1}]\right)\right)^l
\cdot
$$ 
$$ 
\cdot
\gamma_m(\beta)
 \left( v- \left( [a_{i-1}, \dots , a_1]+[a_{i+1}, \dots , a_t]\right)\right)^m \ll
$$ 
$$ 
\ll
\sum_{v=w}^{\infty}
\sum_{\nu+\eta=v}
\sum_{a_1+ \dots +a_t=\nu}
\frac{1}{q \left( a\right)^{2\beta}}
\sum_{a_1+ \dots +a_t'=\eta}
\frac{1}
{ \left( q \left( a\right)q_{-} \left( a\right)\right)^{\beta}}\ll
$$ 
$$ 
\ll
\int_{w}^{\infty}\frac{d v}{v^{\beta-k-1}}=O \left( \frac{1}{w^{\beta-k-2}}\right). 
$$ 
Thus, coefficient at $ k $th term for $ k < \beta-2 $ is equal to
$$ 
R_{k}^{-}=
B_k^{-}+O \left( \frac{1}{w^{\beta-k-2}}\right), 
$$ 
where 
$$ 
B_k^{-}=\sum_{v=1}^{\infty}
\sum_{a \in A_n, a_1+ \dots +a_{i-1}+a_{i+1}+ \dots +a_t=v}
\frac{1}{\langle a_1, \dots , a_{i-1}\rangle^{2\beta}\langle a_{i+1}, \dots , a_t\rangle^{\beta}
\langle a_{i+1}, \dots , a_{t-1}\rangle^{\beta}}\cdot
$$ 
$$ 
\cdot
\sum_{l+m=k}
\gamma_l(\beta)
\left(v- \left( [a_{i-1}, \dots , a_1]+[a_{i+1}, \dots , a_{t-1}]\right)\right)^l
\cdot
\gamma_m(\beta)
\left(v- \left( [a_{i-1}, \dots , a_1]+[a_{i+1}, \dots , a_t]\right)\right)^m. 
$$ 
Now let's consider the terms when $ k \ge \beta-2 $. 
Let's get the estimation for $ R_k^{-} $. 
$$ 
R_k^{-}
\leq
\sum_{v=1}^{w-1}
v^k
\sum_{l+m=k}
\gamma_l(\beta)
\gamma_m(\beta)
$$ 
$$ 
\sum_{\begin{array}{c}
a \in A_n, \\
a_1+ \dots +a_{i-1}+a_{i+1}+ \dots +a_t=v
\end{array}}
\frac{1}{\langle a_1, \dots , a_{i-1}\rangle^{2\beta}\langle a_{i+1}, \dots , a_t\rangle^{\beta}
\langle a_{i+1}, \dots , a_{t-1}\rangle^{\beta}}
\leq
$$ 
$$ 
\leq
\sum_{v=1}^{w-1}
v^k
\sum_{l+m=k}
\gamma_l(\beta)
\gamma_m(\beta)
\sum_{s+u=v}
8C_0\frac{\zeta(2\beta-1)}{\zeta(2\beta)}
\frac{1}{s^{2\beta}u^{\beta}}
\leq
$$ 
$$ 
\leq
8C_0\frac{\zeta(2\beta-1)}{\zeta(2\beta)}
\sum_{l+m=k}
\gamma_l(\beta)
\gamma_m(\beta)
w^{k-\beta+2}. 
$$ 
Then for residual series we get the following estimation:
$$ 
\sum_{k>\beta-2}
\frac{R_{k}^{-}}{n^{2\beta+k}}
\leq
\sum_{k>\beta-2}
\frac{1}{n^{2\beta+k}}
8C_0\frac{\zeta(2\beta-1)}{\zeta(2\beta)}
\sum_{l+m=k}
\gamma_l(\beta)
\gamma_m(\beta)
w^{k-\beta+2}=
$$ 
$$ 
=
8C_0\frac{\zeta(2\beta-1)}{\zeta(2\beta)}
\frac{1}{n^{2\beta}w^{\beta-2}}
\sum_{k>\beta-2}
\frac{w^k}{n^{k}}
\sum_{l+m=k}
\gamma_l(\beta)
\gamma_m(\beta)\leq
$$ 
$$ 
\leq
8C_0\frac{\zeta(2\beta-1)}{\zeta(2\beta)}
\frac{1}{n^{2\beta}w^{\beta-2}}
\sum_{k=1}^{\infty}
\frac{w^k}{n^{k}}
\sum_{l+m=k}
\gamma_l(\beta)
\gamma_m(\beta)=
$$ 
$$ 
=
8C_0\frac{\zeta(2\beta-1)}{\zeta(2\beta)}
\frac{1}{n^{2\beta}w^{\beta-2}}
\left(
\frac{1}
{1-
\frac{w}{n}
}
\right)^{2\beta}. 
$$ 
With the given $ w $ magnitude $ \frac{1}{1-\frac{w}{n}} $ 
 doesn't exceed 2, thus , residual series can be estimated as 
 $ O(\frac{1}{n^{2\beta}w^{\beta-2}}) $. 
If $ \beta $ is integer, then for $ k=\beta-2 $ we obtain
$$ 
\sum_{v=1}^{w-1}
\sum_{a \in A_n, a_1+ \dots +a_{i-1}+a_{i+1}+ \dots +a_t=v}
\frac{1}{\langle a_1, \dots , a_{i-1}\rangle^{2\beta}\langle a_{i+1}, \dots , a_t\rangle^{\beta}
\langle a_{i+1}, \dots , a_{t-1}\rangle^{\beta}}
$$ 
$$ 
\sum_{l+m=k}
 \gamma_l(\beta)
 \left( v- \left( [a_{i-1}, \dots , a_1]+[a_{i+1}, \dots , a_{t-1}]\right)\right)^l
\cdot
$$ 
$$ 
\cdot
 \gamma_m(\beta)
 \left( v- \left( [a_{i-1}, \dots , a_1]+[a_{i+1}, \dots , a_t]\right)\right)^m
=
O \left( \log w\right). 
$$ 
Similar actions can be made for part of the sum with $ \frac{1}{qq_{+}} $. 
\\We get
$$ 
\sum_{v=1}^{w-1}
\sum_{a \in A_n, a_1+ \dots +a_{i-1}+a_{i+1}+ \dots +a_t=v}
\frac{1}{ \left( qq_{+}\right)^{\beta}}=
$$ 
$$ 
=\frac{1}{n^{2\beta}}
\sum_{v=1}^{w-1}
\sum_{a \in A_n, a_1+ \dots +a_{i-1}+a_{i+1}+ \dots +a_t=v}
\frac{1}{\langle a_1, \dots , a_{i-1}\rangle^{2\beta}\langle a_{i+1}, \dots , a_t\rangle^{\beta}
\langle a_{i+1}, \dots , a_t-1\rangle^{\beta}}
+
$$ 
$$ 
\sum_{k=1}^{\infty}
\frac{1}{n^{2\beta+k}}
\sum_{v=1}^{w-1}
\sum_{a \in A_n, a_1+ \dots +a_{i-1}+a_{i+1}+ \dots +a_t=v}
\frac{1}{\langle a_1, \dots , a_{i-1}\rangle^{2\beta}\langle a_{i+1}, \dots , a_t\rangle^{\beta}
\langle a_{i+1}, \dots , a_t-1\rangle^{\beta}}
$$ 
$$ 
\sum_{l+m=k}
 \gamma_l(\beta)
 \left( v- \left( [a_{i-1}, \dots , a_1]+[a_{i+1}, \dots , a_t-1]\right)\right)^l
\cdot
 \gamma_m(\beta)
 \left( v- \left( [a_{i-1}, \dots , a_1]+[a_{i+1}, \dots , a_t]\right)\right)^m. 
$$ 
Coefficient $ R_{k}^{+} $ at $ k $th term equals
$$ 
R_{k}^{+}=
(\sum_{v=1}^{\infty}-\sum_{v=w}^{\infty})
\sum_{a \in A_n, a_1+ \dots +a_{i-1}+a_{i+1}+ \dots +a_t=v}
\frac{1}{\langle a_1, \dots , a_{i-1}\rangle^{2\beta}\langle a_{i+1}, \dots , a_t\rangle^{\beta}
\langle a_{i+1}, \dots , a_t-1\rangle^{\beta}}
$$ 
$$ 
\sum_{l+m=k}
 \gamma_l(\beta)
 \left( v- \left( [a_{i-1}, \dots , a_1]+[a_{i+1}, \dots , a_t-1]\right)\right)^l
\cdot
 \gamma_m(\beta)
 \left( v- \left( [a_{i-1}, \dots , a_1]+[a_{i+1}, \dots , a_t]\right)\right)^m. 
$$ 
Let us consider the first series. 
It can be majorized with th following series:
$$ 
\sum_{v=1}^{\infty}v^k
\sum_{\nu+\eta=v}
\sum_{a_1+ \dots +a_t=\nu}
\frac{1}{q \left( a\right)^{2\beta}}
\sum_{a_1+ \dots +a_t'=\eta}
\frac{1}
{ \left( q \left( a\right)q_{+} \left( a\right)\right)^{\beta}}. 
$$ 
According to lemma 9, 
$$ 
\sum_{a_1+ \dots +a_t'=\nu}
\frac{1}
{ \left( q \left( a\right)\right)^{2\beta}}=O \left( \frac{1}{\nu^{2\beta}}\right), 
$$ 
and it follows from lemma 6 in [3] that
$$ 
\sum_{a_1+ \dots +a_t'=\eta}
\frac{1}
{ \left( q \left( a\right)q_{+} \left( a\right)\right)^{\beta}}=O \left( \frac{1}{\eta^{2\beta}}\right), 
$$ 
then common term of this series can be estimated as $ O \left( \frac{1}{v^{2\beta-k-1}}\right) $, 
so, the series converges when $ k<2\beta-2 $. \\
With this $ k $ let us estimate the residual series:
$$ 
\sum_{v=w}^{\infty}
\sum_{a \in A_n, a_1+ \dots +a_{i-1}+a_{i+1}+ \dots +a_t=v}
\frac{1}{\langle a_1, \dots , a_{i-1}\rangle^{2\beta}\langle a_{i+1}, \dots , a_t\rangle^{\beta}
\langle a_{i+1}, \dots , a_t-1\rangle^{\beta}}
$$ 
$$ 
\sum_{l+m=k}
 \gamma_l(\beta)
 \left( v- \left( [a_{i-1}, \dots , a_1]+[a_{i+1}, \dots , a_t-1]\right)\right)^l
\cdot
$$ 
$$ 
\cdot
\gamma_m(\beta)
 \left( v- \left( [a_{i-1}, \dots , a_1]+[a_{i+1}, \dots , a_t]\right)\right)^m \ll
$$ 
$$ 
\ll
\sum_{v=w}^{\infty}
\sum_{\nu+\eta=v}
\sum_{a_1+ \dots +a_t=\nu}
\frac{1}{q \left( a\right)^{2\beta}}
\sum_{a_1+ \dots +a_t'=\eta}
\frac{1}
{ \left( q \left( a\right)q_{+} \left( a\right)\right)^{\beta}}\ll
$$ 
$$ 
\ll
\int_{w}^{\infty}\frac{d v}{v^{2\beta-k-1}}=O \left( \frac{1}{w^{2\beta-k-2}}\right). 
$$ 
Thus, coefficient $ R_{k}^{+} $ at $ k $th term when $ k < 2\beta-2 $ equals
$$ 
R_{k}^{+}=
B_k^{+}+O \left( \frac{1}{w^{2\beta-k-2}}\right), 
$$ 
where 
$$ 
B_k^{+}=\sum_{v=1}^{\infty}
\sum_{a \in A_n, a_1+ \dots +a_{i-1}+a_{i+1}+ \dots +a_t=v}
\frac{1}{\langle a_1, \dots , a_{i-1}\rangle^{2\beta}\langle a_{i+1}, \dots , a_t\rangle^{\beta}
\langle a_{i+1}, \dots , a_t-1\rangle^{\beta}}
$$ 
$$ 
\sum_{l+m=k}
\gamma_l(\beta)
\left(v- \left( [a_{i-1}, \dots , a_1]+[a_{i+1}, \dots , a_t-1]\right)\right)^l
\cdot
\gamma_m(\beta)
\left(v- \left( [a_{i-1}, \dots , a_1]+[a_{i+1}, \dots , a_t]\right)\right)^m. 
$$ 
Now let's consider the terms when $ k \ge 2\beta-2 $. 
Let's get the estimation for $ k $th term . 
$$ 
R_k^{+}
\leq
\sum_{v=1}^{w-1}
v^k
\sum_{l+m=k}
\gamma_l(\beta)
\gamma_m(\beta)
$$ 
$$ 
\sum_{\begin{array}{c}
a \in A_n, \\
a_1+ \dots +a_{i-1}+a_{i+1}+ \dots +a_t=v
\end{array}
}
\frac{1}{\langle a_1, \dots , a_{i-1}\rangle^{2\beta}\langle a_{i+1}, \dots , a_t\rangle^{\beta}
\langle a_{i+1}, \dots , a_t-1\rangle^{\beta}}
\leq
$$ 
$$ 
\leq
\sum_{v=1}^{w-1}
v^k
\sum_{l+m=k}
\gamma_l(\beta)
\gamma_m(\beta)
\sum_{s+u=v}
4C_0
\frac{1}{s^{2\beta}u^{2\beta}}
\leq
$$ 
$$ 
\leq
4C_0
\sum_{l+m=k}
\gamma_l(\beta)
\gamma_m(\beta)
w^{k-2\beta+2}. 
$$ 
Then we get the following estimation for the residual series. 
$$ 
\sum_{k>2\beta-2}
\frac{R_{k}^{+}}{n^{2\beta+k}}
\leq
\sum_{k>2\beta-2}
\frac{1}{n^{2\beta+k}}
4\left(
\frac{\zeta(2\beta-1)}{\zeta(2\beta)}+
\left(\frac{\zeta(2\beta-1)}{\zeta(2\beta)}\right)^2
\right)
\sum_{l+m=k}
\gamma_l(\beta)
\gamma_m(\beta)
w^{k-2\beta+2}=
$$ 
$$ 
=
4C_0
\frac{1}{n^{2\beta}w^{2\beta-2}}
\sum_{k>2\beta-2}
\frac{w^k}{n^{k}}
\sum_{l+m=k}
\gamma_l(\beta)
\gamma_m(\beta)\leq
$$ 
$$ 
\leq
4C_0
\frac{1}{n^{2\beta}w^{2\beta-2}}
\sum_{k=1}^{\infty}
\frac{w^k}{n^{k}}
\sum_{l+m=k}
\gamma_l(\beta)
\gamma_m(\beta)
\leq
$$ 
$$ 
\leq
4C_0
\frac{1}{n^{2\beta}w^{\beta-2}}
\left(
\frac{1}
{1-
\frac{w}{n}
}
\right)^{2\beta}. 
$$ 
With the given $ w $ magnitude $ \frac{1}
{1-
\frac{w}{n}
} $ doesn't exceed 2, so, residual series can be estimated as 
 $ O(\frac{1}{n^{2\beta}w^{2\beta-2}}) $. 
If $ 2\beta $ is integer, then for $ k=2\beta-2 $ we get
$$ 
\sum_{v=1}^{w-1}
\sum_{\begin{array}{c}
a \in A_n, \\
a_1+ \dots +a_{i-1}+a_{i+1}+ \dots +a_t=v
\end{array}}
\frac{1}{\langle a_1, \dots , a_{i-1}\rangle^{2\beta}\langle a_{i+1}, \dots , a_t\rangle^{\beta}
\langle a_{i+1}, \dots , a_t-1\rangle^{\beta}}
$$ 
$$ 
\sum_{l+m=k}
 \gamma_l(\beta)
 \left( v- \left( [a_{i-1}, \dots , a_1]+[a_{i+1}, \dots , a_t-1]\right)\right)^l
\cdot
$$ 
$$ 
\cdot
 \gamma_m(\beta)
 \left( v- \left( [a_{i-1}, \dots , a_1]+[a_{i+1}, \dots , a_t]\right)\right)^m
=
O \left( \log w\right). 
$$ 
\\Thus, 
adding sum for $ \frac{1}{q \left( a\right)q_{-} \left( a\right)} $ to sum for $ \frac{1}{q \left( a\right)q_{+} \left( a\right)} $, 
we obtain when $ \beta $ is integer 
$$ 
\Sigma_{n, 2}^{ \left( 3\right)}
=
\sum_{0\le k<\beta-2}
B_k
\frac{1}{n^{2\beta+k}}
+O \left( 
\frac{1}{n^{ \left( \beta-1\right) \left( 2r-1\right)}}
+
\frac{\log w}{n^{3\beta-2}}
+
\frac{1}{w^{\beta-2}n^{2\beta}}
\right), 
$$ 
when $ \beta $ is not integer
$$ 
\Sigma_{n, 2}^{ \left( 3\right)}
=
\sum_{0\le k<\beta-2}
B_k
\frac{1}{n^{2\beta+k}}
+O \left( 
\frac{1}{n^{ \left( \beta-1\right) \left( 2r-1\right)}}
+
\frac{1}{w^{\beta-2}n^{2\beta}}
\right), 
$$ 
where 
$$ 
B_k=B_k^{-}+B_k^{+}. 
$$ 
{\bf Lemma is proved. }

{\bf Lemma 14. }
{\it
When $ 2\beta $ is integer 
\begin{equation}
\label{sigma_n13+_int}
\Sigma_{n, 1}^{ \left( 3 \right)+}
=
\frac{1}{n^{2\beta}}\frac{2\zeta \left( 2\beta-1 \right)}{\zeta \left( 2\beta \right)}
+
\sum_{1\le k<2\beta-1}D_k\frac{1}{n^{2\beta+k}}
+O \left( 
\frac{1}{n^{2\beta}w^{2 \left( \beta-1 \right)}}
+
\frac{\log w}{n^{4\beta-1}}
+
\frac{1}{n^{ \left( \beta-1 \right) \left( 2r-1 \right)}}
+
\frac{1}{w^{2\beta-1}n^{2\beta}}
 \right), 
\end{equation}
when $ 2\beta $ is not integer 
\begin{equation}
\label{sigma_n13+}
\Sigma_{n, 1}^{ \left( 3 \right)+}
=
\frac{1}{n^{2\beta}}\frac{2\zeta \left( 2\beta-1 \right)}{\zeta \left( 2\beta \right)}
+
\sum_{1\le k<2\beta-1}D_k\frac{1}{n^{2\beta+k}}
+O \left( 
\frac{1}{n^{2\beta}w^{2 \left( \beta-1 \right)}}
+
\frac{1}{n^{ \left( \beta-1 \right) \left( 2r-1 \right)}}
+
\frac{1}{w^{2\beta-1}n^{2\beta}}
 \right). 
\end{equation}
where $ D_k $ are some constants. }

{\bf Proof. }
\newline
$$ 
\Sigma_{n, 1}^{ \left( 3 \right)+}
=
\sum_{a \in A_n, q \left( a \right) < n^r, a_t > n-w}
\frac{1}{ \left( qq_{+} \right)^{\beta}}
$$ 
$$ 
=
\sum_{a \in A_n, a_t > n-w}
\frac{1}{ \left( qq_{+} \right)^{\beta}}
-
\sum_{a \in A_n, q \left( a \right) \geq n^r}
\frac{1}{ \left( qq_{+} \right)^{\beta}}
$$ 
The second sum can be estimated according to lemma 11 as 
$$ 
\sum_{a \in A_n, q \left( a \right) \geq n^r}
\frac{1}{ \left( qq_{+} \right)^{\beta}}=O \left( 
\frac{1}{n^{ \left( \beta-1 \right) \left( 2r-1 \right)}}
 \right). 
$$ 
For the first sum we have
\begin{equation}
\label{sum+}
\sum_{a \in A_n, a_t > n-w}
\frac{1}{ \left( qq_{+} \right)^{\beta}}
=
\sum_{v=1}^{w}
\sum_{a \in A_n, a_1 + \ldots + a_{t-1}=v}
\frac{1}{ \left( qq_{+} \right)^{\beta}}. 
\end{equation}
Here $ a_t=n-v, v=1, . . , w $. 
As 
$$ 
q = a_t q_{-} + \left( q_{-} \right)_{-} 
=
q_{-} \left( n-v+[a_{t-1}, \dots , a_1] \right), 
$$ 
$$ 
q_{+} = \left( a_t-1 \right)q_{-} + \left( q_{-} \right)_{-} 
=
q_{-} \left( n-v -1+[a_{t-1}, \dots , a_1] \right), 
$$ 
then we get for $ \frac{1}{ \left( qq_{+} \right)^{\beta}} $ 
\begin{equation}
\label{undersum}
\frac{1}{ \left( qq_{+} \right)^{\beta}}
=
\frac{1}{n^{2\beta}q_{-}^{2\beta}}
 \left( 
\frac{1}{ \left( 1-\frac{v}{n}+\frac{1}{n}[a_{t-1}, \dots , a_1] \right)^{\beta}}\cdot
\frac{1}{ \left( 1-\frac{v}{n}-\frac{1}{n}+\frac{1}{n}[a_{t-1}, \dots , a_1] \right)^{\beta}}
 \right). 
\end{equation}
Expanding magnitudes
 $ 
\frac{1}{ \left( 1-\frac{v}{n}+\frac{1}{n}[a_{t-1}, \dots , a_1] \right)^{\beta}}
 $ 
and
 $ 
\frac{1}{ \left( 1-\frac{v}{n}-\frac{1}{n}+\frac{1}{n}[a_{t-1}, \dots , a_1] \right)^{\beta}}
 $ 
into Teilor series according to parameters $ \frac{v}{n}-\frac{1}{n}[a_{t-1}, \dots , a_1] $ we obtain
 $ \frac{v}{n}+\frac{1}{n}-\frac{1}{n}[a_{t-1}, \dots , a_1] $ correspondingly , 
when $ v \le n-1 $ 
$$ 
\frac{1}{ \left( 1-\frac{v}{n}+\frac{1}{n}[a_{t-1}, \dots , a_1] \right)^{\beta}}=
1+
\sum_{k=1}^{\infty}
 \gamma_{k}(\beta)
 \left( \frac{v}{n}-\frac{1}{n}[a_{t-1}, \dots , a_1] \right)^k, 
$$ 
$$ 
\frac{1}{ \left( 1-\frac{v}{n}-\frac{1}{n}+\frac{1}{n}[a_{t-1}, \dots , a_1] \right)^{\beta}}=
1+
\sum_{k=1}^{\infty}
 \gamma_{k}(\beta)
 \left( \frac{v}{n}+\frac{1}{n}-\frac{1}{n}[a_{t-1}, \dots , a_1] \right)^k, 
$$ 
where $ \gamma_k(\beta) $ are defined in (\ref{gamma}). 
Then, substituting obtained series into (\ref{undersum}), 
we get when $ v \le n-1 $ 
$$ 
\frac{1}{ \left( qq_{+} \right)^{\beta}}
=
$$ 
$$ 
=
\frac{1}{n^{2\beta}q_{-}^{2\beta}}
 \left( 1 + 
\sum_{k=1}^{\infty}
\sum_{l+m=k}
 \gamma_l(\beta)
 \left( \frac{v}{n}-\frac{1}{n}[a_{t-1}, \dots , a_1] \right)^l
\cdot
\gamma_m(\beta)
 \left( \frac{v}{n}+\frac{1}{n}-\frac{1}{n}[a_{t-1}, \dots , a_1] \right)^m
 \right)=
$$ 
$$ 
=\frac{1}{n^{2\beta}q_{-}^{2\beta}}
 \left( 1 + 
\sum_{k=1}^{\infty}\frac{1}{n^{k}}
\sum_{l+m=k}
 \gamma_l(\beta)
 \left( v-[a_{t-1}, \dots , a_1] \right)^l
\cdot
\gamma_m(\beta)
 \left( v+1-[a_{t-1}, \dots , a_1] \right)^m
 \right).
$$ 
Next, substituting expression for $ \frac{1}{ \left( qq_{+} \right)^{\beta}} $ into (\ref{sum+}), we obtain
$$ 
\sum_{v=1}^{w-1}\sum_{a_1+ \dots +a_{t-1}=v}\frac{1}{ \left( qq_{+} \right)^{\beta}}
=
\sum_{v=1}^{w-1}\sum_{a_1+ \dots +a_{t-1}=v}
\frac{1}{n^{2\beta}q_{-}^{2\beta}}
\cdot
$$
$$
\cdot
 \left( 1 + 
\sum_{k=1}^{\infty}\frac{1}{n^{k}}
\sum_{l+m=k}
 \gamma_l(\beta)
 \left( v-[a_{t-1}, \dots , a_1] \right)^l
\cdot
\gamma_m(\beta)
 \left( v+1-[a_{t-1}, \dots , a_1] \right)^m
 \right)=
$$ 
$$ 
=
\frac{2}{n^{2\beta}}
\sum_{v=1}^{w-1}\sum_{a_1+ \dots +a_{t-1}=v, a_{t-1} \geq 2}
\frac{1}{q_{-}^{2\beta}}
 + 
%$$  
%$$  
 %+
\sum_{k=1}^{\infty}
\frac{2}{n^{2\beta+k}}
\sum_{v=1}^{w-1}\sum_{a_1+ \dots +a_{t-1}=v, a_{t-1} \geq 2}
\frac{1}{q_{-}^{2\beta}} \cdot
$$
$$
\cdot
\sum_{l+m=k}
 \gamma_l(\beta)
 \left( v-[a_{t-1}, \dots , a_1] \right)^l
\cdot
\gamma_m(\beta)
 \left( v+1-[a_{t-1}, \dots , a_1] \right)^m. 
$$ 
Let's consider sum $ R_k $ at $ \frac{1}{n^{2\beta+k}} $. 
$$ 
R_k=\sum_{v=1}^{w-1}\sum_{a_1+ \dots +a_{t-1}=v, a_{t-1} \geq 2}
\frac{2}{q_{-}^{2\beta}}
\sum_{l+m=k}
 \gamma_l(\beta)
 \left( v-[a_{t-1}, \dots , a_1] \right)^l
\cdot
\gamma_m(\beta)
 \left( v+1-[a_{t-1}, \dots , a_1] \right)^m 
=
$$ 
$$ 
= \left( \sum_{v=1}^{\infty}-\sum_{v=w}^{\infty} \right)\sum_{a_1+ \dots +a_{t-1}=v, a_{t-1} \geq 2}
\frac{2}{q_{-}^{2\beta}}
\sum_{l+m=k}
 \gamma_l(\beta)
 \left( v-[a_{t-1}, \dots , a_1] \right)^l
\cdot
\gamma_m(\beta)
 \left( v+1-[a_{t-1}, \dots , a_1] \right)^m. 
$$ 
Let us consider the first sum:
$$ 
\sum_{v=1}^{\infty}\sum_{a_1+ \dots +a_{t-1}=v, a_{t-1} \geq 2}
\frac{2}{q_{-}^{2\beta}}
\sum_{l+m=k}
 \gamma_l(\beta)
 \left( v-[a_{t-1}, \dots , a_1] \right)^l
%\cdot
%$$  
%$$  
\cdot
\gamma_m(\beta)
 \left( v+1-[a_{t-1}, \dots , a_1] \right)^m \ll 
$$ 
$$ 
\ll
\sum_{v=1}^{\infty}v^k
\sum_{a_1+ \dots +a_{t-1}=v, a_{t-1} \geq 2}
\frac{1}{q_{-}^{2\beta}}.
$$ 
According to lemma 9
$$ 
\sum_{a_1+ \dots +a_t=\nu}
\frac{1}{q \left( a \right)^{2\beta}}=O \left( \frac{1}{\nu^{2\beta}} \right), 
$$ 
thus, common term of the given series is $ O \left( \frac{1}{n^{2\beta-k}} \right) $, so
series converges when $ 2\beta-k>1 $, the given when $ k<2\beta-1 $. 
Let us estimate residual series with these $ k $ :
$$ 
\sum_{v=w}^{\infty}\sum_{a_1+ \dots +a_{t-1}=v}
\frac{2}{q_{-}^{2\beta}}
\sum_{l+m=k}
 \gamma_l(\beta)
 \left( v-[a_{t-1}, \dots , a_1] \right)^l
\cdot
\gamma_m(\beta)
 \left( v+1-[a_{t-1}, \dots , a_1] \right)^m \ll 
$$ 
$$ 
\ll
\sum_{v=w}^{\infty}v^k
\sum_{a_1+ \dots +a_i=v}
\frac{1}{q^{2\beta}}
\ll
\int_{w}^{\infty}\frac{d v}{v^{2\beta-k}}=O \left( \frac{1}{w^{2\beta-k-1}} \right). 
$$ 
Thus, 
$$ 
R_k=D_k+O \left( \frac{1}{w^{2\beta-k-1}} \right), 
$$ 
where 
$$ 
D_k=\sum_{v=1}^{\infty}\sum_{a_1+ \dots +a_{t-1}=v, a_{t-1} \geq 2}
\frac{2}{q_{-}^{2\beta}}
\sum_{l+m=k}
 \gamma_l(\beta)
 \left( v-[a_{t-1}, \dots , a_1] \right)^l
\cdot
\gamma_m(\beta)
 \left( v+1-[a_{t-1}, \dots , a_1] \right)^m. 
$$ 
Now let us consider $ R_k $ when $ k \ge 2\beta-1 $ :
$$ 
R_k
\leq
2\sum_{l+m=k}
 \gamma_l(\beta)\cdot
\gamma_m(\beta)
\sum_{v=1}^{w-1}v^k
\sum_{a_1+ \dots +a_i=v, a_{t-1} \geq 2}
\frac{1}{q^{2\beta}}\leq
$$ 
$$ 
\leq
4\frac{\zeta(2\beta-1)}{\zeta(2\beta)}
\sum_{l+m=k}
 \gamma_l(\beta)\cdot
\gamma_m(\beta)
\int_{1}^{w-1}\frac{d v}{v^{2\beta-k}}\leq
4\frac{\zeta(2\beta-1)}{\zeta(2\beta)}
\sum_{l+m=k}
 \gamma_l(\beta)\cdot
\gamma_m(\beta)
w^{k+1-2\beta}, 
$$ 
if $ k>2\beta-1 $, and $ O \left( \log w \right) $ when $ k = 2\beta-1 $ (in the case when $ 2\beta $ is integer ). 
Then for residual series we get the following estimation:
$$ 
\sum_{k>2\beta-1}
\frac{R_k}{n^{2\beta+k}}
\leq
4\frac{\zeta(2\beta-1)}{\zeta(2\beta)}
\sum_{k>2\beta-1}
\frac{w^{k+1-2\beta}}{n^{2\beta+k}}
\sum_{l+m=k}
 \gamma_l(\beta)\cdot
\gamma_m(\beta)=
$$ 
$$ 
=
4\frac{\zeta(2\beta-1)}{\zeta(2\beta)}
\frac{1}{w^{2\beta-1}n^{2\beta}}
\sum_{k>2\beta-1}
\frac{w^{k}}{n^{k}}
\sum_{l+m=k}
 \gamma_l(\beta)\cdot
\gamma_m(\beta)\leq
$$ 
$$ 
\leq
4\frac{\zeta(2\beta-1)}{\zeta(2\beta)}
\frac{1}{w^{2\beta-1}n^{2\beta}}
\sum_{k=1}^{\infty}
\frac{w^{k}}{n^{k}}
\sum_{l+m=k}
 \gamma_l(\beta)\cdot
\gamma_m(\beta)\leq
$$ 
$$ 
\leq
4\frac{\zeta(2\beta-1)}{\zeta(2\beta)}
\frac{1}{w^{2\beta-1}n^{2\beta}}
\left(\frac{1}{1-\frac{w}{n}}\right)^{2\beta}
=O\left(\frac{1}{w^{2\beta-1}n^{2\beta}}
\right)
$$ 
with the given $ w $. \\
Extracting the constant in the main term , we obtain when $ 2\beta $ is integer
$$ 
\Sigma_{n, 1}^{ \left( 3 \right)+}
=
\frac{1}{n^{2\beta}}\frac{2\zeta \left( 2\beta-1 \right)}{\zeta \left( 2\beta \right)}
+
\sum_{1\le k<2\beta-1}D_k\frac{1}{n^{2\beta+k}}
+O \left( 
\frac{1}{n^{2\beta}w^{2 \left( \beta-1 \right)}}
+
\frac{\log w}{n^{4\beta-1}}
+
\frac{1}{n^{ \left( \beta-1 \right) \left( 2r-1 \right)}}
+
\frac{1}{w^{2\beta-1}n^{2\beta}}
 \right), 
$$ 
when $ 2\beta $ is not integer 
$$ 
\Sigma_{n, 1}^{ \left( 3 \right)+}
=
\frac{1}{n^{2\beta}}\frac{2\zeta \left( 2\beta-1 \right)}{\zeta \left( 2\beta \right)}
+
\sum_{1\le k<2\beta-1}D_k\frac{1}{n^{2\beta+k}}
+O \left( 
\frac{1}{n^{2\beta}w^{2 \left( \beta-1 \right)}}
+
\frac{1}{n^{ \left( \beta-1 \right) \left( 2r-1 \right)}}
+
\frac{1}{w^{2\beta-1}n^{2\beta}}
 \right). 
$$ 
{\bf Lemma is proved. }
\\
{\bf Lemma 15. }
{\it
When $ 2\beta $ is integer
\begin{equation}
\label{main_term_integer}
\Sigma_{n, 1}^{ \left( 3 \right)-}
=
\frac{1}{n^{\beta}}
\frac{2\zeta \left( 2\beta-1 \right)}{\zeta \left( 2\beta \right)}
+\sum_{1\le k< 2\beta-1}
E_k\frac{1}{n^{\beta+k}}
+O \left( 
\frac{\log w}{n^{3\beta-1}}
+
\frac{1}{n^{\beta}w^{2 \left( \beta-1 \right)}}
 \right), 
\end{equation}
otherwise
\begin{equation}
\label{main_term_nonint}
\Sigma_{n, 1}^{ \left( 3 \right)-}
=
\frac{1}{n^{\beta}}
\frac{2\zeta \left( 2\beta-1 \right)}{\zeta \left( 2\beta \right)}
+\sum_{1\le k< 2\beta-1}
E_k\frac{1}{n^{\beta+k}}
+O \left( 
\frac{1}{n^{\beta}w^{2 \left( \beta-1 \right)}}
 \right), 
\end{equation}
where $ E_k $ are some constants. 
}\\
{\bf Proof. }
$$ 
\Sigma_{n, 1}^{ \left( 3 \right)-}
=
\sum_{a \in A_n, q \left( a \right) < n^r, a_t > n-w}
\frac{1}{ \left( qq_{-} \right)^{\beta}}
$$ 
$$ 
=
\sum_{a \in A_n, a_t > n-w}
\frac{1}{ \left( qq_{-} \right)^{\beta}}
-
\sum_{a \in A_n, q \left( a \right) \geq n^r}
\frac{1}{ \left( qq_{-} \right)^{\beta}}
$$ 
The second sum can be estimated according to lemma 11 as 
$$ 
\sum_{a \in A_n, q \left( a \right) \geq n^r}
\frac{1}{ \left( qq_{-} \right)^{\beta}}=O \left( 
\frac{1}{n^{ \left( \beta-1 \right) \left( 2r-1 \right)}}
 \right). 
$$ 
For the first sum we have
\begin{equation}
\label{main_sum}
\sum_{a \in A_n, a_t > n-w}
\frac{1}{ \left( qq_{-} \right)^{\beta}}
=
2\sum_{v=1}^{w}
\sum_{
\begin{array}{c}
a \in A_n, \\
a_1 + \ldots + a_{t-1}=v, \\
a_{t-1}\geq 2
\end{array}
}
\frac{1}{ \left( qq_{-} \right)^{\beta}}. 
\end{equation}
Here $ a_t=n-v, v=1, . . , w $. 
As 
$$ 
q = a_t q_{-} + \left( q_{-} \right)_{-} 
=
nq_{-} \left( 1-\frac{1}{n}(v-[a_{t-1}, \dots , a_1]) \right), 
$$ 
then expanding $ \frac{1}{(qq_{-})^{\beta}} $ into Teilor series according to parameter
 $ \frac{1}{n}(v-[a_{t-1}, \dots , a_1]) $, we get 
\begin{equation}
\label{teilor}
\frac{1}{ \left( qq_{-} \right)^{\beta}}
=
\frac{1}{n^{\beta}q_{-}^{2\beta}}
+
\frac{1}{n^{\beta}q_{-}^{2\beta}}
\sum_{k=1}^{\infty}
\gamma_k(\beta)\left(\frac{1}{n}(v-[a_{t-1}, \dots , a_1])\right)^k. 
\end{equation}
Thus, substituting (\ref{teilor}) into (\ref{main_sum}), we obtain
$$ 
\sum_{a \in A_n, a_t > n-w}
\frac{1}{ \left( qq_{-} \right)^{\beta}}
=
2\sum_{v=1}^{w}
\sum_{
\begin{array}{c}
a \in A_n, \\
a_1 + \ldots + a_{t-1}=v, \\
a_{t-1}\geq 2
\end{array}
}
\frac{1}{n^{\beta}q_{-}^{2\beta}}
+
$$
$$
+
\sum_{k=1}^{\infty}
\frac{1}{n^{\beta+k}}
2\gamma_k(\beta)
\sum_{v=1}^{w}
\left(v-[a_{t-1}, \dots , a_1]\right)^k
%\cdot
%$$  
%$$  
%\cdot
\sum_{
\begin{array}{c}
a \in A_n, \\
a_1 + \ldots + a_{t-1}=v, \\
a_{t-1}\geq 2
\end{array}
}
\frac{1}{q_{-}^{2\beta}}
= 
$$ 
$$
=
2\sum_{v=1}^{w}
\sum_{
\begin{array}{c}
a \in A_n, \\
a_1 + \ldots + a_{t-1}=v, \\
a_{t-1}\geq 2
\end{array}
}
\frac{1}{n^{\beta}q_{-}^{2\beta}}
+
\sum_{k=1}^{\infty}
\frac{R_k}{n^{\beta+k}},
$$
%Let's consider sum at $ \frac{1}{n^{\beta+k}} $. 
where $R_k$ is defined as follows:
$$ 
R_k= 2\gamma_k(\beta)
\sum_{v=1}^{w}
\sum_{a \in A_v}
\frac{1}{q^{2\beta}}
 \left( v-[a_t, \dots , a_1] \right)^{k}= 
$$ 
$$ 
=
 2\gamma_k(\beta)
 \left( \sum_{v=1}^{\infty}-\sum_{v=w}^{\infty} \right)
\sum_{a \in A_v}
\frac{1}{q^{2\beta}}
 \left( v-[a_t, \dots , a_1]\right)^{k}. 
$$ 
Let us consider the first sum:
$$ 
 2\gamma_k(\beta)
 \sum_{v=1}^{\infty}
\sum_{a \in A_v}
\frac{1}{q^{2\beta}}
 \left( v-[a_t, \dots , a_1] \right)^{k} \ll 
%$$  
%$$  
%\ll
\sum_{v=1}^{\infty}
v^{k} 
\sum_{a \in A_v}
\frac{1}{q^{2\beta}}. 
$$ 
It follows from lemma 9 that 
 $ \sum_{a \in A_v}
\frac{1}{q^{2\beta}}=O \left( \frac{1}{v^{2\beta}} \right)
 $, hence 
$$ 
\sum_{v=1}^{\infty}v^k
\sum_{a \in A_v}
\frac{1}{q^{2\beta}} \ll
\sum_{v=1}^{\infty}\frac{1}{v^{2\beta-k}}. 
$$ 
Thus, the given series converges when $ 2\beta-k>1 $, i. e. when $ k<2\beta-1 $. 
With these $ k $ let us estimate the residual series of the given series. 
$$ 
 2\gamma_k(\beta)
\sum_{v=w}^{\infty}
\sum_{a \in A_v}
\frac{1}{q^{2\beta}}
 \left( v-[a_t, \dots , a_1] \right)^{k} \ll 
$$ 
$$ 
\ll
\sum_{v=w}^{\infty}
v^{k} 
\sum_{a \in A_v}
\frac{1}{q^{2\beta}} \ll 
\int_{w}^{\infty}\frac{d v}{v^{2\beta-k}}=O \left( \frac{1}{w^{2\beta-k-1}} \right). 
$$ 
Hence, when $ k<2\beta-1 $
$$
R_k=E_k+O\left(\frac{1}{w^{2\beta-k-1}}\right).
$$
Here $E_k$ are constants, defined with the following formula
$$
E_k=
2\gamma_k(\beta)
 \sum_{v=1}^{\infty}
\sum_{a \in A_v}
\frac{1}{q^{2\beta}}
 \left( v-[a_t, \dots , a_1] \right)^{k}.
$$
Now let's estimate the sum when $ k \ge 2\beta-1 $. 
$$ 
R_k 
\leq
4\gamma_k(\beta)
C_0
\sum_{v=1}^{w-1}
v^{k-2\beta} 
\leq 
%\int_{1}^{w-1}\frac{d v}{v^{2\beta-k}}
$$ 
$$ 
\leq
4\gamma_k(\beta)
C_0
w^{k+1-2\beta} 
$$ 
when $ k>2\beta-1 $ and $ O \left( \log w \right) $ when $ k=2\beta-1 $. 
Then, summing according to $ k > 2\beta-1 $, we get
$$ 
\sum_{k > 2\beta-1}
\frac{R_k}{n^{\beta+k}} \leq
4
C_0
\frac{1}{n^{\beta}w^{2\beta-1}}
\sum_{k > 2\beta-1}
\gamma_k(\beta)
\frac{w^k}{n^k} \leq
$$ 
$$ 
\leq
4
C_0
\frac{1}{n^{\beta}w^{2\beta-1}}
\left(
\frac{1}{
1-\frac{w}{n} 
}
\right)^{\beta}=
O\left(
\frac{1}{n^{\beta}w^{2\beta-1}}
\right). 
$$ 
\\Extracting the constant in the main term , we get
$$ 
\frac{2}{n^\beta}
\sum_{a \in A_n, a_1 + \ldots + a_{t-1} \leq w}
\frac{1}{q_{-}^{2\beta}} = 
\frac{1}{n^\beta}
\frac{2\zeta \left( 2\beta-1 \right)}{\zeta \left( 2\beta \right)}+
O \left( \frac{1}{n^{\beta}w^{2 \left( \beta-1 \right)}} \right). 
$$ 
Thus, when $ 2\beta $ is integer
$$ 
\Sigma_{n, 1}^{ \left( 3 \right)-}
=
\frac{1}{n^{\beta}}
\frac{2\zeta \left( 2\beta-1 \right)}{\zeta \left( 2\beta \right)}
+\sum_{1\le k<2\beta-1}E_k\frac{1}{n^{\beta+k}}
+O \left( 
\frac{\log n}{n^{3\beta-1}}
+
\frac{1}{n^{\beta}w^{2 \left( \beta-1 \right)}}
 \right), 
$$ 
when $ 2\beta $is not integer 
$$ 
\Sigma_{n, 1}^{ \left( 3 \right)-}
=
\frac{1}{n^{\beta}}
\frac{2\zeta \left( 2\beta-1 \right)}{\zeta \left( 2\beta \right)}
+\sum_{1\le k<2\beta-1}E_k\frac{1}{n^{\beta+k}}
+O \left( 
\frac{1}{n^{\beta}w^{2 \left( \beta-1 \right)}}
 \right). 
$$ 

{\bf Lemma is proved. }

{\bf The final step in proving theorem 2. }
Substituting (\ref{big_q}), (\ref{sigma_n22}), (\ref{sigma_n23_int}), (\ref{sigma_n23_int}), 
(\ref{sigma_n13+_int}), (\ref{sigma_n13+_int}), (\ref{main_term_integer}), (\ref{main_term_nonint})
into (\ref{whole_sum}), 
we obtain when $ 2\beta $ is integer 
$$ 
\Sigma_{n, 1}^{ \left( 3 \right)-}
+\Sigma_{n, 1}^{ \left( 3 \right)+}
+\Sigma_{n, 2}^{ \left( 3 \right)}
+\Sigma_{n, 2}^{ \left( 2 \right)}
+\Sigma_{n, 2}^{ \left( 1 \right)}
=
$$ 
$$ 
=
\frac{1}{n^{\beta}}
\frac{2\zeta \left( 2\beta-1 \right)}{\zeta \left( 2\beta \right)}
+\sum_{1\le k< 2\beta-2}{C}_k\frac{1}{n^{\beta+k}}+
+\sum_{0\le k< \beta-2}{C^*}_k\frac{1}{n^{2\beta+k}}+
$$ 
$$ 
+O \left( 
\frac{\log n}{n^{3\beta-2}}
+
\frac{1}{n^{2\beta}w^{\beta-2}}
+
\frac{ n^2\left( \log^{3\beta} n \right)^{3\beta}}{w^{3\beta}}
+
\frac{1}{n^{\beta}w^{2 \left( \beta-1 \right)}}
%\frac{1}{n^{\beta}w^{2 \left( \beta-1 \right)}}
%+\frac{\log n}{n^{3\beta-1}}
+
\frac{1}{n^{(\beta-1)(2r-1)}}
 \right), 
$$ 
where $ C_k=E_k $, 
and $ C^{*}_k=B_{k}+D_k$,
$ {C^*}_0=B_k+\frac{2\zeta(2\beta-1)}{\zeta(2\beta)}$.  
When $2 \beta $ is not integer , 
$$ 
\Sigma_{n, 1}^{ \left( 3 \right)-}
+\Sigma_{n, 1}^{ \left( 3 \right)+}
+\Sigma_{n, 2}^{ \left( 3 \right)}
+\Sigma_{n, 2}^{ \left( 2 \right)}
+\Sigma_{n, 2}^{ \left( 1 \right)}
=
$$
$$
=
\frac{1}{n^{\beta}}
\frac{2\zeta \left( 2\beta-1 \right)}{\zeta \left( 2\beta \right)}
+\sum_{1\le k< 2\beta-2}{C}_k\frac{1}{n^{\beta+k}}
+\sum_{0\le k< \beta-2}{C^*}_k\frac{1}{n^{2\beta+k}}+
$$ 
$$ 
+O \left( 
\frac{1}{n^{2\beta}w^{\beta-2}}
+
\frac{ n^2\left( \log^{3\beta} n \right)^{3\beta}}{w^{3\beta}}
+
\frac{1}{n^{\beta}w^{2 \left( \beta-1 \right)}}
+
\frac{1}{n^{(\beta-1)(2r-1)}}
 \right), 
$$ 
where $ C_k=E_k $, and $ {C^*}_k=B_k+D_k $, 
$ {C^*}_0=B_k+\frac{2\zeta(2\beta-1)}{\zeta(2\beta)}$. 
To minimize degree of error term, let $ w=\frac{n}{2}-2 $, 
 $ r=\frac{3\beta}{2(\beta-1)}+\frac{1}{2} $. \\
{\bf Theorem is proved. }\\

\end{document}